\documentclass[]{amsart}
\usepackage[utf8]{inputenc}
\usepackage[mathscr]{eucal} % So-called Euler fonts, allows \mathscr
\usepackage{stmaryrd} % St. Mary's font, fixes a few bad symbols
\usepackage{amsmath,amsthm,amsfonts,amssymb,amscd} % Standard AMS packages
\usepackage[usenames,dvipsnames]{xcolor} % Allows colors
\usepackage{xspace} % Allows smart spacing after period
\usepackage{fancyhdr} % For footers, headers
\usepackage{graphicx} % For inserting images
\usepackage{listings} % For source code listing
\usepackage[]{hyperref} % For hyperlinks
\usepackage{enumitem} % More control for enumerated lists
\usepackage{relsize} % Change font size relative to current size
\usepackage{tikz-cd} % For commutative diagrams
\usepackage{mdwlist} % Additional list customization
\usepackage{multicol} % Allows multiple columns
\usepackage{float} % Improved placement for floating objects
\usepackage{adjustbox} % Adjust boxed content
\usepackage{tikz} % Regular tikz
\usepackage[noabbrev,nameinlink]{cleveref} % Additional reference customization
\usepackage{setspace} % Used for changing spacing
\usepackage{bbm} % more bb
\usetikzlibrary{matrix, calc, arrows}

\hypersetup{
	colorlinks=true, %set true if you want colored links
	linktoc=all,     %set to all if you want both sections and subsections linked
	linkcolor=blue,  %choose some color if you want links to stand out
}

\setlist[enumerate]{label=(\alph*)} % Default enumerate alpha
\usetikzlibrary{graphs,decorations.pathmorphing,decorations.markings}
\tikzcdset{scale cd/.style={every label/.append style={scale=#1},
		cells={nodes={scale=#1}}}}

\newcommand{\Hom}{\operatorname{Hom}}
\newcommand{\Aut}{\operatorname{Aut}}

\newcommand{\inv}{^{-1}}

\newcommand{\op}{^{\operatorname{op}}}

\newcommand{\id}{\operatorname{id}}

\newcommand{\im}{\operatorname{im}}

\newcommand{\stab}{\operatorname{stab}}

\newcommand{\Br}{\operatorname{Br}}

\newcommand{\Spc}{\operatorname{Spc}}

\newcommand{\cone}{\operatorname{cone}}

\newcommand{\on}[1]{\operatorname{#1}}

\newcommand{\calI}{\mathcal{I}}

\newcommand{\calP}{\mathcal{P}}

\newcommand{\calX}{\mathcal{X}}

\newcommand{\catB}{\mathscr{B}}
\newcommand{\catC}{\mathscr{C}}
\newcommand{\catD}{\mathscr{D}}

\newcommand{\catI}{\mathscr{I}}
\newcommand{\catJ}{\mathscr{J}}
\newcommand{\catK}{\mathscr{K}}
\newcommand{\catL}{\mathscr{L}}
\newcommand{\catM}{\mathscr{M}}

\newcommand{\catP}{\mathscr{P}}
\newcommand{\catQ}{\mathscr{Q}}
\newcommand{\catR}{\mathscr{R}}
\newcommand{\catS}{\mathscr{S}}

\newcommand{\frakp}{\mathfrak{p}}

\newcommand{\bbone}{\mathbbm{1}}

\newcommand{\ev}{\on{ev}}
\newcommand{\coev}{\on{coev}}
\newcommand{\Sp}{\on{Sp}}

\newcommand{\supp}{\on{supp}}
\newcommand{\sbull}{{\scriptscriptstyle\bullet}}

\newcommand{\Ext}{\on{Ext}}
\newcommand{\Spec}{\on{Spec}}
\newcommand{\cp}{^{cp}}

\newcommand{\Proj}{\on{Proj}}

\newtheorem{theorem}{Theorem}[section]
\newtheorem{lemma}[theorem]{Lemma}
\newtheorem{prop}[theorem]{Proposition}
\newtheorem{corollary}[theorem]{Corollary}

\newtheorem*{theorem*}{Theorem}

\theoremstyle{remark}
\newtheorem{remark}[theorem]{Remark}

\theoremstyle{definition}
\newtheorem{definition}[theorem]{Definition}

\makeatletter
\newcommand{\xhookdoubleheadrightarrow}[2][]{%
	\lhook\joinrel
	\ext@arrow 0359\rightarrowfill@ {#1}{#2}%
	\mathrel{\mspace{-15mu}}\rightarrow
}
\makeatother

\newtheorem{innercustomthm}{Theorem}
\newenvironment{customthm}[1]
{\renewcommand\theinnercustomthm{#1}\innercustomthm}
{\endinnercustomthm}

\begin{document}
	\title{On functoriality and the tensor product property in noncommutative tensor-triangular geometry}
	\author{Sam K. Miller}
	\address{Department of Mathematics, University of California Santa Cruz, Santa Cruz, CA 95064} %required
	% \email{sakmille@ucsc.edu} %optional
	 %required
	 %optional
	\begin{abstract}
		Two pertinent questions for any support theory of a monoidal triangulated category are whether it is functorial and if the tensor product property holds. To this end, we consider the complete prime spectrum of an essentially small monoidal triangulated category, which we show is universal among support data satisfying the tensor product property, even if it is empty. The complete prime spectrum is functorial and parametrizes radical thick tensor-ideals, a noncommutative analogue of Balmer's reconstruction theorem. We give criteria for when induced maps on complete prime spectra are injective or surjective, and determine the complete prime spectrum for crossed product categories. Finally, we determine the universal functorial support theory for monoidal triangulated categories coinciding with the Balmer spectrum on braided monoidal triangulated categories.   
	\end{abstract}
	
	\maketitle
	\markright{\MakeUppercase{On functoriality and the tensor product property in nc-ttg}}
	
	\section{Introduction}
	
	Support is a fundamental notion in algebraic geometry, representation theory, and far beyond. Classically, the support of a module $M$ over a commutative ring $R$, denoted $\supp(M)$, is the subspace of $\Spec(R)$ of all prime ideals $\frakp$ for which the localization $M_\frakp \neq 0$. Generalizations of this notion developed over time, first appearing in modular representation theory via the work of Quillen \cite{Q71}, Carlson \cite{Car83}, and many others, towards understanding cohomological properties of modules. A key construction in the modular representation theoretic setting is the \textit{cohomological support variety} of a $kG$-module $M$, a subvariety $V_G(M)$ of the projective variety $V_G = \Proj(\on{H}^\sbull(G,k))$. These varieties communicate cohomological data such as periodicity or projectivity of a module. Moreover, Benson, Carlson, and Rickard \cite{BCR97} demonstrated that the classification of thick tensor-ideals of the tensor-triangulated category $\stab(kG)$, the stable module category for $kG$, could be read off from $V_G$. 
	
	Support theories and classifications such as these were united by Balmer via tensor-triangular geometry. Balmer constructed what is now referred to as the Balmer spectrum $\Spc(\catK)$ of a ``tensor-triangulated'' (by this, we mean a symmetric monoidal triangulated) category $\catK$, the set of all prime thick $\otimes$-ideals of $\catK$ \cite{Bal05}. Each object $x \in \catK$ comes equipped with a support $\supp(x) \subset \catK$, an abstract ``support variety'' for the object $x$, giving $\Spc(\catK)$ the structure of a locally ringed space. Balmer's key insight was that one can read off all the thick tensor-ideals of $\catK$ from the pair $(\Spc(\catK), \supp)$, akin to Benson-Carlson-Rickard's theorem for the cohomological support variety $V_G$. Moreover, the Balmer spectrum is the \textit{universal support} for essentially small categories - the pair $(\Spc(\catK), \supp)$ is the final object in the category of support theories satisfying the following properties:
	
	\begin{enumerate}
		\item $\supp(0) = \emptyset$ and $\supp(\bbone) = \Spc(\catK)$;
		\item $\supp(a \oplus b) = \supp(a) \cup \supp(b)$;
		\item $\supp(\Sigma a) = \supp(a)$, where $\Sigma$ denotes the translation of $\catK$;
		\item $\supp(a) \subseteq \supp(b) \cup \supp(c)$ for any exact triangle $a \to b \to c \to \Sigma a$;
		\item $\supp(a \otimes b) = \supp(a) \cap \supp(b)$.
	\end{enumerate}
	
	These properties generalize the classical properties of the support of a module over a commutative ring, and frequently hold for any well-studied support variety in a tensor-triangular setting, such as the cohomological support varieties $V_G(M)$ of $kG$-modules $M$ (whether this holds for any cohomological support variety of a braided tensor category is posed in \cite{BPW24}). Circling back to modular representation theory, the cohomological support variety $V_G$ is recovered as the Balmer spectrum of the stable module category $\stab(kG)$. Tensor-triangular geometry has since evolved to a highly active field of (meta)mathematics, with implications in representation theory, stable homotopy theory, motivic geometry, and beyond. For a more in-depth overview of tensor-triangular geometry, we point the reader to \cite{Bal10b}. 
	
	But what about the oft-occurring monoidal triangulated categories for which the tensor product is not necessarily symmetric? Analogously to the differences between standard algebraic geometry vs. noncommutative algebraic geometry (see \cite{Rey25} for a reader-friendly overview), things are less clear-cut. Support varieties for settings without a symmetric tensor product (or any tensor product at all) have been studied for quite some time \cite{EHSST03, FW11, GS23, Lin10, SS04} and Buan, Krause, and Solberg laid some of the framework for developing noncommutative tensor-triangular geometry via theories of support using frame theory in \cite{BKS07, BKSS20}. However, only recently a noncommutative analogue of the Balmer spectrum was formalized by Nakano, Vashaw, and Yakimov in \cite{NVY22}, and further examined by the authors in \cite{NVY22b, NVY23, NVY24}. A similar construction was adopted by Negron and Pevtsova in \cite{NP23}. 
	
	Nakano, Vashaw, and Yakimov follow the blueprint laid out by both Buan-Krause-Solberg and noncommutative algebraic geometry, by defining the noncommutative Balmer spectrum $\Spc(\catK)$ to be the set of all \textit{prime thick two-sided tensor-ideals}, i.e. thick tensor ideals $\catP$ for which $\catI \otimes \catJ \subseteq \catP$ implies $\catI \subseteq \catP$ or $\catJ \subseteq \catP$. They construct the supports of objects of $\catK$ analogously, and from there, draw up axioms for noncommutative support data and weak support data, and prove the noncommutative Balmer spectrum is the universal support. The noncommutative support data is defined via a weaker tensor product axiom than the standard counterpart: for noncommutative support data, we require the weaker condition \[\supp(a) \cap \supp(b) = \bigcup_{c\in \catK}\supp(a \otimes c\otimes b),\] which is equivalent to the condition $\supp(a \otimes b) = \supp(a) \cap \supp(b)$ when $\catK$ is symmetric, or more generally, braided. The authors also lay out a strategy for the classification of thick one-sided and two-sided tensor ideals via classifying supports, in a similar manner as the symmetric setting. Analysis of the Balmer spectrum is also amenable towards using frame-theoretic methods, as considered in \cite{BKS07, BKS20, MR23, GS23} - in fact, the Balmer spectrum frame-theoretically appears to be the correct choice for classifying thick $\otimes$-ideals. All seems to be smooth sailing. 
	
	However, in noncommutative algebraic geometry, there is not just one noncommutative spectrum - there are numerous, all of which are imperfect in the sense of Reyes \cite{Rey12}, in that either a spectral theory is non-functorial, or it cannot be guaranteed to be empty. In the noncommutative tensor-triangular setting, one might expect a similar snafu. The noncommutative Balmer spectrum $\Spc(\catK)$ as a topological space does not have many of the well-behaved properties that its commutative counterpart has - for instance if $\Spc(\catK)$ is non-Noetherian, it is not known if closed points, i.e. minimal primes, exist or if $\Spc(\catK)$ is quasi-compact. Perhaps more critically, $\Spc$ is in general non-functorial, which eliminates many of the techniques we have for tensor-triangular classification in the noncommutative setting! This confirms the expectation \cite[Remark 3.3.3]{NVY22}.
	
	\begin{customthm}{A}[Theorem \ref{ex:nocps}]
		Let $\catK = D_b(k) \oplus D_b(k)$. The crossed product category $\catK \rtimes C_2$, with $C_2$ acting by permuting indices, contains no completely prime ideals. Moreover, the inclusion $F\colon  \catK \hookrightarrow \catK \rtimes C_2$ does not induce a well-defined map $\Spc(F): \Spc(\catK \rtimes C_2) \to \Spc(\catK)$ nor a well-defined functor (see Definition \ref{def:catofsupp}) $\catS(F): \catS(\catK \rtimes C_2) \to \catS(\catK)$. 
	\end{customthm}
	
	More generally, we give a criteria for precisely when a monoidal exact functor $F\colon  \catK \to \catL$ induces a well-defined map on Balmer spectra. As a result, we can deduce the noncommutative Balmer spectra of localizations and idempotent completion, adaptations of standard tensor-triangular geometric results for the noncommutative setting. 
	
	\begin{customthm}{B}[Theorem \ref{cor:wheniswelldef}]
		Let $F\colon \catK \to\catL$ be an exact $\otimes$-functor, and denote its image by $F(\catK)\subseteq \catL$. Then $F$ induces a continuous map on Balmer spectra $\varphi := \Spc(F): \Spc(\catL) \to \Spc(\catK)$ given by $\catQ \mapsto F\inv(\catQ)$ if and only if for all prime thick $\otimes$-ideals $\catQ \in \Spc(\catL)$, the intersection $F(\catK) \cap \catQ$ satisfies the primality property, \[\text{For all $a, b \in F(K)$, } a \otimes F(K) \otimes c \in F(K) \cap \catQ \text{ implies } a \in \catQ \text{ or } b\in \catQ.\]
	\end{customthm}   
	
	The lack of functoriality is perhaps discouraging, but there is more to the picture. First, Nakano, Vashaw, and Yakimov propose alternative axioms of support to consider, encouraging the notion that there may not be one fixed notion of support data, akin to how there are multiple spectral theories for noncommutative rings, or how support axioms break down for big tensor-triangulated categories. In their companion piece \cite{NVY22b}, the authors analyze when $(\Spc(\catK), \supp)$ satisfies the tensor product property, i.e. $\supp(a) \cap \supp(b) = \supp(a \otimes b)$. They prove this property holds if and only if every prime is \textit{completely prime}, i.e. a prime $\catP$ which satisfies $a \otimes b \in \catP$ implies $a \in \catP$ or $b \in \catP$. The question of whether the the tensor product property holds for cohomological support varieties is of independent interest and has been examined further in \cite{NP23, BPW24, BPW25}.
	
	If we think of this complete prime spectrum $\Spc\cp(\catK)$ as another construction of a spectral space or support theory, then in some ways it is better-behaved; for instance, it is indeed functorial. However, it has its own major shortcomings (as explained in \cite{NVY22}) - notably, it may be empty, and the standard ``prime lifting'' techniques (see \cite[Lemma 2.2]{Bal05} or \cite[Theorem 3.2.3]{NVY22}) do not seem feasible if $\Spc\cp(\catK) \subset \Spc(\catK)$. However, in tensor-triangular geometry, one should think about the Balmer spectrum in the terms of universal support: the (commutative) Balmer spectrum is the universal support datum for tensor-triangulated categories, and the noncommutative Balmer spectrum is universal for, in fact, multiple notions of noncommutative support data. And similarly, the complete prime spectrum is itself universal - it is the final object in the category of \textit{multiplicative support data}, i.e. support data satisfying the tensor product property. 
	
	\begin{customthm}{C}
		Let $\catS^\otimes(\catK)$ denote the category of multiplicative support data (see Definition \ref{def:supportdata} and \ref{def:catofsupp}). The following hold.
		
		\begin{itemize}
			\item $\catS^\otimes(-): \mathbf{mon}_\Delta\op \to \mathbf{Cat}$ is a contravariant functor from the category of essentially small monoidal triangulated categories to the category of essentially small categories; (Theorem \ref{thm:functorial}) 
			\item  The complete prime spectrum $(\Spc\cp(\catK), \supp\cp)$ is the final object of $\catS^\otimes(\catK)$, even if $\Spc\cp(\catK) = \emptyset$; (Theorem \ref{thm:thm6})
			\item The assignment $\Spc\cp(-): \mathbf{mon}_\Delta\op \to \mathbf{Top}$ is a contravariant functor from the category of essentially small monoidal triangulated categories to the category of topological spaces. (Theorem \ref{thm:thm5}) 
		\end{itemize}
		
	\end{customthm}
	
	Given a spectral theory on a monoidal triangulated category $\catK$, the specialization-closed or Thomason subsets should be expected to parametrize some data about $\catK$. Reframing slightly, we prove an analogue of Balmer's reconstruction theorem \cite[Theorem 4.10]{Bal05}. As in the commutative setting, the complete primes parametrize the radical thick $\otimes$-ideals of $\catK$. In fact, a complete prime is nothing more than a radical prime. 
	
	\begin{customthm}{D}
		Let $\Spc\cp(\catK)^\vee$ denote the Hochster dual of $\Spc\cp(\catK)$ (see Definition \ref{def:hochster}). 
		\begin{enumerate}
			\item Let $\catP$ be a thick $\otimes$-ideal of $\catK$. Then $\catP$ is completely prime if and only if $\catP$ is prime and radical; (Proposition \ref{prop:cpprimeradical})
			\item We have an order-preserving bijection from the set $\catR$ of radical thick $\otimes$-ideals to the set $\catS$ of \textit{open} subsets of $\Spc\cp(\catK)^\vee$ induced by the assignment \[\catI \mapsto \bigcup_{x \in \catI} \supp(x).\] In particular, every radical thick $\otimes$-ideal of $\catK$ is an intersection of completely prime thick $\otimes$-ideals of $\catK$. (Theorem \ref{thm:cpparametrization})
		\end{enumerate}
	\end{customthm}
	
	We determine the complete prime spectrum for the crossed product category $\catK \rtimes G$, using the classification of $\Spc(\catK \rtimes G)$ by Huang and Vashaw \cite{HV25}. In addition, we compute the one-sided thick $\otimes$-ideals of $\catK \rtimes G$. 
	
	\begin{customthm}{E}[Theorem \ref{thm:crossedcprimes}]
		Let $G$ be a group acting on $\catK$, that is, we have a functor $G \to \Aut_{\otimes, \Delta}(\catK)$. Then the homeomorphism $\Spc(\catK \rtimes G) \cong G\on{-}\Spc(\catK)$ (see \cite[Proposition 5.4]{HV25}) restricts to a homeomorphism $\Spc\cp(\catK \rtimes G) \cong G\on{-}\Spc\cp(\catK)$. 
	\end{customthm}
	
	We adapt a result of Balmer \cite{Bal18} to the noncommutative setting regarding induced maps on spectra $\Spc\cp(F): \Spc\cp(\catL) \to \Spc\cp(\catK)$ under the additional assumption that $\catL$ is \textit{duo}, i.e. right, left, and two-sided thick tensor-ideals coincide. 
	
	\begin{customthm}{F}[Theorem \ref{thm:thm2}]
		Suppose $\catK$ is rigid and $\catL$ is a duo monoidal triangulated category. Suppose that the exact $\otimes$-functor $F\colon \catK \to \catL$ detects $\otimes$-nilpotents of morphisms, i.e. every $f\colon x\to y$ in $\catK$ such that $F(f) = 0$ satisfies $f^{\otimes n} = 0$ for some $n \geq 1$. Then the induced map $\varphi\colon \Spc\cp(\catL)\to \Spc\cp(\catK)$ is surjective.
	\end{customthm}
	
	Finally, following the blueprint laid out by Reyes \cite{Rey12}, we expand our scope and determine the universal functorial support theory for monoidal triangulated categories that coincides with $\Spc$ for tensor-triangulated categories, the partial prime spectrum. Moreover, the partial prime spectrum universal in the sense of support axioms - it is the final \textit{commutative support datum} (see Definition \ref{def:commsupp}). 
	
	\begin{customthm}{G}
		We say a full subcategory $\catP$ of $\catK$ is a \textit{partial prime $\otimes$-ideal} if for all (not necessarily thick or full) monoidal triangulated subcategories $\catC$ of $\catK$ which can be given a braided structure, $\catP \cap \catC$ is a prime thick $\otimes$-ideal of $\catC$. Let $p\Spc(\catK)$ denote the set all partial prime $\otimes$-ideals of $\catK$.
		\begin{itemize}
			\item  The functor $p\Spc: \mathbf{mon}_\Delta\op \to \mathbf{Top}$ is the final object in the full subcategory of $\on{Fun}(\mathbf{mon}_\Delta\op, \mathbf{Top})$ consisting of functors $F\colon  \mathbf{mon}_\Delta\op \to\mathbf{Top}$ which are naturally isomorphic to $\Spc$ upon restriction to $\mathbf{Bmon}_\Delta\op$; (Corollary \ref{cor:universalfunctorialfunctor})
			\item  If $p\Spc(\catK) \neq \emptyset$, then the pair $(p\Spc(\catK), p\supp)$ is a final object in the category $\catS^{cm}(\catK)$ of commutative support data (see Definition \ref{def:commsupp}). (Theorem \ref{thm:finalcommsuppdata})
		\end{itemize}
	\end{customthm}
	
	The paper is organized as follows. Section 2 discusses preliminaries in noncommutative tensor-triangular geometry and proves a few topological properties analogous to the symmetric setting. Section 3 introduces the category of support data, proves its partial functoriality, and contains discussion about multiplicative support axioms. Section 4 gives a criterion for when monoidal exact functors induce maps on noncommutative Balmer spectrum, and deduces a few standard yet important properties about the spectrum with regards to localization and idempotent completion. Section 5 proves functoriality and universality of the complete prime spectrum. Section 6 proves an analogue to Balmer's classification of thick tensor-ideals for complete primes. Section 7 is a brief aside: we deduce the final object in the category of quasi-support data, as defined in \cite{NVY22}. 
	
	In Section 8, we apply the support theories and Balmer spectra we have considered to the crossed product category $\catK \rtimes G$, an extension of Benson and Witherspoon's crossed product algebra \cite{BW14}. In the next two sections, we switch to analyzing the functoriality of the complete prime spectrum, with the additional assumption that one of our categories is \textit{duo}, i.e. every one-sided ideal is two-sided. In Section 9, we introduce this notion, discuss rigidity, and discuss the work of Negron-Pevtsova \cite{NP23}. In Section 10, we adapt some of Balmer's surjectivity results \cite{Bal18} to the noncommutative setting for the complete prime spectrum. Finally in Section 11, we define the partial prime spectrum and determine that it is the universal functorial support theory for monoidal triangulated categories, analogous to Reyes \cite{Rey12}.

	\textbf{Acknowledgments:} The author thanks Julia Pevtsova for numerous conversations which led to certain key aspects of the paper, Kent Vashaw, Dan Nakano, and Timothy de Deyn for discussions on aspects of noncommutative tensor-triangular geometry, and Beren Sanders for his detailed feedback on an earlier version of this paper. He also is indebted to Paul Balmer and Martin Gallauer, whose correspondence and conversations throughout the past year helped him quickly enter the world of tensor-triangular geometry. 
	
	\section{Preliminaries}
	
	We begin by reviewing standard noncommutative tensor-triangular preliminaries, noncommutative analogues of tensor-triangular geometry developed in \cite{Bal05}. We follow the conventions of \cite{NVY22b}. 
	
	\begin{definition}\label{def:primes}
		A \textit{monoidal triangulated category} (M$\Delta$C for short) $\catK$ is a monoidal category in the sense of \cite[Definition 2.2.1]{EGNO15} which is triangulated and for which the monoidal structure $\otimes\colon  \catK \times \catK \to \catK$ is an exact bifunctor. For this paper, we assume $\catK$ and $\catL$ are essentially small monoidal triangulated categories unless we state otherwise. 
		
		A \textit{thick subcategory} of a triangulated category $\catK$ is a full triangulated subcategory of $\catK$ that contains all direct summands. Note that it follows that any thick subcategory of $\catK$ is \textit{replete}, i.e. closed under isomorphism classes. A thick \textit{right} (resp. \textit{left}, resp. \textit{two-sided}) \textit{tensor ideal} ($\otimes$-ideal for short) of $\catK$ is a thick triangulated subcategory of $\catK$ that is closed under right tensoring (resp. left, resp. left and right tensoring) with arbitrary objects of $\catK$. For any $x \in \catK$, we may speak of \textit{the right (resp. left, resp. two-sided) ideal generated by $x$}, this is the minimal right (resp. left, resp. two-sided) ideal containing $x$. These will be denoted by $\langle x\rangle_r$, $\langle x \rangle_l$, and $\langle x\rangle$ respectively. We omit the ``two-sided'' when we speak of two-sided ideals for brevity - unless specified, ideals will be assumed two-sided. 
		
		A proper thick $\otimes$-ideal $\catP \subset \catK$ is \textit{prime} if \[\catI \otimes \catJ \subseteq \catP \implies \catI \subseteq \catP \text{ or } \catJ \subseteq \catP\] for all ideals $\catI, \catJ$. In fact, this property is equivalent to saying that this holds for all pairs of right ideals $\catI, \catJ$. It is also equivalent to the condition that for all $x,y \in \catK$, \[x \otimes z \otimes y \in \catP \text{ for all } z \in \catK \implies x \in \catP \text{ or } y \in \catP.\]
		A thick $\otimes$-ideal $\catP$ is \textit{completely prime} if it satisfies the ``usual'' element-wise primality condition, \[x \otimes y \in \catP \implies x \in \catP \text{ or } y \in \catP.\]        
	\end{definition}

	\begin{definition}\label{def:spc}
		The \textit{noncommutative Balmer spectrum} $\Spc(\catK)$ is the set of thick prime $\otimes$-ideals. For every object $x \in \catK$, we define its \textit{support} by \[\supp(x) := \{\catP \in \Spc(\catK) \mid x \not\in \catP\}.\] The \textit{Zariski topology} on $\Spc(\catK)$ is the one generated by the following open subsets $\{U(x) := \Spc(\catK) \setminus \supp(x) \mid x \in \catK\}$. 
		
		We let $\Spc^{cp}(\catK)$ be the topological subspace consisting of all completely prime ideals of $\catK$, the \textit{noncommutative complete prime spectrum}, or \textit{CP-spectrum} for short. Set \[\supp^{cp}(x) = \supp(x) \cap \Spc^{cp}(\catK),\] then the Zariski topology on $\Spc^{cp}(\catK)$ is the one generated by compliments $U^{cp}(x)$ of $\supp^{cp}(x)$ over all $x \in \catK$. We omit the ``Balmer'' from the title to avoid wordiness and confusion. 
		
		Both $\Spc(\catK)$ and $\Spc\cp(\catK)$ have a Hochster-dual topology, with topology generated by the open sets $\supp(x)$ and $\supp\cp(x)$ respectively for all $x \in \catK$. We say any subset of $\Spc(\catK)$ or $\Spc\cp(\catK)$ is \textit{dual-open} if it is open in the Hochster-dual topology.
	\end{definition}
	
	\begin{prop}\label{prop:t0}
		For any point $\catP \in \Spc^i(\catK)$, with $i\in \{\emptyset, cp\}$, its closure is \[\overline{\{\catP\}} = \{\catQ \in \Spc^i(\catK) \mid \catQ \subseteq \catP\}.\] In particular, if $\overline{\{\catP_1\}} = \overline{\{\catP_2\}}$, then $\catP_1 = \catP_2$, i.e. the space $\Spc^i(\catK)$ is $T_0$. 
	\end{prop}
	\begin{proof}
		Let $S_0 = \catK \setminus \catP$. We set $Z(S_0) := \{\catP \in \Spc^i(\catK) \mid S_0 \cap \catP = \emptyset\}$ (in particular, $Z(\{x\}) = \supp^i(x)$). Then $\catP \in Z(S_0)$, and if $\catP \in Z(S)$ for some $S \subset \catK$, then $S \subseteq S_0$, hence $Z(S_0) \subseteq Z(S)$. Therefore $Z(S_0)$ is the smallest subset containing $\catP$. But it is easy to see $Z(S_0) = \bigcap_{x\in S_0} \supp^i(x)$, so $Z(S_0)$ is closed. Thus, $Z(S_0)$ is the smallest closed subset containing $\catP$, i.e. $\{\catP\} = Z(S_0) = \{\catQ \in \Spc^i(\catK)\mid \catQ \subset \catP\}$. The final statement is clear. 
	\end{proof}

	\begin{prop}{\cite[Theorem 3.2.3]{NVY22}}\label{prop:primelift}
		Given a two-sided $\otimes$-multiplicative class $S$ of objects in $\catK$ and a proper thick two-sided $\otimes$-ideal $\catI \subset \catK$ such that $\catI \cap S = \emptyset$, then any maximal element of the set \[X(S, \catI) := \{\catJ \text{ a thick $\otimes$-ideal of $\catK$} \mid \catI \subseteq \catJ, \catJ \cap S = \emptyset\}\] is a prime ideal. In particular, there exists a prime $\catP \in \Spc(\catK)$ such that $\catI \subseteq \catP$ and $\catP \cap S = \emptyset$, and $\Spc(\catK)$ is nonempty. 
	\end{prop}
	
	The following results of prime lifting are easily shown.
	
	\begin{prop}{\cite[Proposition 2.3]{Bal05}}
		\begin{enumerate}
			\item Let $S$ be a two-sided $\otimes$-multiplicative collection of objects which does not contain zero. There exists a prime ideal $\catP \in \Spc(\catK)$ such that $\catP \cap S = \emptyset$.
			\item Let $\catI \subset \catK$ be a proper thick $\otimes$-ideal. Then there exists a maximal proper thick $\otimes$-ideal $\catM \subset \catK$ containing $\catI$.
			\item Maximal proper thick $\otimes$-ideals are prime. 
			\item The spectrum $\Spc(\catK)$ is nonempty. 
		\end{enumerate}
	\end{prop}
	
	\begin{corollary}{\cite[Corollary 2.4, 2.5]{Bal05}}
		An object $x \in \catK$ belongs to all primes, i.e. $\supp(x) = \emptyset$, if and only if it is $\otimes$-nilpotent, and belongs to no prime, i.e. $\supp(x) = \Spc(\catK)$, if and only if it generates $\catK$ as a thick two-sided $\otimes$-ideal. 
	\end{corollary}

	Analogous to the classical setting in \cite{Bal05}, the pair $(\Spc(\catK), \supp)$ is universal in the following sense, as was shown in \cite{NVY22}. Given a topological space $Y$, let $\calX(Y)$, $\calX_{cl}(Y)$, and $\calX_{sp}(Y)$ denote the collections of its subsets, closed subsets, and specialization closed subsets respectively. 
	
	\begin{definition}\label{def:supportdata}
		A \textit{support datum} (as defined in \cite{NVY22}) for a monoidal triangulated category $\catK$ is a map \[\sigma\colon \catK \to \calX_{cl}(Y)\] for a topological space $Y$ such that the following hold:
		\begin{enumerate}
			\item $\sigma(0) = \emptyset$ and $\sigma(\bbone) = Y$;
			\item $\sigma(x \oplus y) = \sigma(x) \cup \sigma(y)$;
			\item $\sigma(\sum x) = \sigma(x)$;
			\item $\sigma(x) \subseteq \sigma(y) \cup \sigma(z)$ for all distinguished triangles $x \to y \to z \to \sum x$;
			\item $\bigcup_{z \in \catK}\sigma(x \otimes z \otimes y) = \sigma(x) \cap \sigma(y)$.
		\end{enumerate}
		A \textit{morphism of support data} is a continuous map $f\colon (Y,\sigma)\to (Z,\tau)$ such that $f\inv(\tau(x)) = \sigma(x)$ for all $x\in \catK$. 
		
		Any support datum satisfying the following stronger condition \[ (e') \,\,\sigma(x \otimes y) = \sigma(x) \cap \sigma(y)\] is said to have the \textit{tensor product property}, and we call such a support data satisfying (e') a \textit{multiplicative support datum}. This property is present in Balmer's support axioms in \cite{Bal05}, and it is easy to verify that for a braided monoidal triangulated category, (e) and (e') are equivalent. 
		
		By default, we allow $(\emptyset, \id)$ to be a multiplicative support datum. 
	\end{definition}
	
	The noncommutative Balmer spectrum is universal, as it should be. 
	
	\begin{theorem}{\cite[Theorem 4.2.2]{NVY22}}
		The tuple $(\Spc(\catK),\supp)$ is the final object in the category of support data $(Y, \sigma)$ for $\catK$ such that $\sigma(x)$ is closed for each $x \in \catK$. For any such $\sigma\colon \catK \to \calX(Y)$ there is a unique continuous map $f\colon Y \to \Spc(\catK)$ satisfying \[\sigma(x) = f\inv(\supp(x)).\]
	\end{theorem}
	
	We will show $\Spc\cp$ is similarly universal, but first, we observe that it is indeed a multiplicative support datum. 
	
	\begin{prop}\label{prop:multsuppdata}
		The pair $(\Spc^{cp}, \supp^{cp})$ is a multiplicative support datum. 
	\end{prop}
	\begin{proof}
		Verifications of (a)-(d) of the support data axioms follows the same as $(\Spc, \supp)$. Verification of (e) follows analogously to the symmetric case, we spell out the details. We verify the dual statement, that $U^{cp}(x \otimes y) = U^{cp}(x) \cup U^{cp}(Y)$. By definition, $U^{cp}(x) = \{\catP \in \Spc^{cp}\mid x \in \catP\}$. Since $\catP$ is a tensor ideal, for every $x \in \catP$, $x \otimes y, y \otimes x \in \catP$, hence $U(x) \cup U(y) \subseteq U(x \otimes y)$. On the other hand, since $\catP$ is completely prime, if $x \otimes y \in\catP$, then either $x \in \catP$ or $y \in \catP$, so $U(x \otimes y) \subseteq U(x) \cup U(y)$, as desired. 
	\end{proof}
	
	Furthermore, from the proof it is easy to see that if there exists a prime $\catP \in \Spc(\catK) \setminus \Spc\cp(\catK)$, then there exist objects $x, y \in \catK$ satisfying $U(x) \otimes U(y) \subset U(x\otimes y)$. Indeed, choose $x, y$ such that $x\otimes y \in \catP$ but $x, y \not\in\catP$. The next results follow essentially the same as the proofs in \cite{Bal05}. Note that $\Spc^{cp}(\catK)$ may be empty (see Theorem \ref{ex:nocps}).  
	
	\begin{prop}{\cite[Proposition 2.11, Corollary 2.12]{Bal05}}\label{prop:closedpts}
		If $\Spc(\catK)$ is nonzero, there exist minimal primes in $\catK$. In particular, any nonempty closed subset of $\Spc^{cp}(\catK)$ contains at least one closed point. 
	\end{prop}
	
	\begin{prop}{\cite[Proposition 2.14(b)]{Bal05}}
		Assume $\Spc\cp(\catK)$ is nonzero. Any quasi-compact open of $\Spc\cp(\catK)$ is of the form $U(x)$ for some $x \in \catK$.
	\end{prop}
	\begin{proof}
		Let $U$ be a quasi-compact open. Then $U \subseteq \bigcup_{x\in \catK} U(x)$ is an open cover of $U$, hence there exists a finite subcover $U = U(x_1) \cup \dots \cup U(x_n) = U(x_1 \otimes \cdots \otimes x_n)$, as desired. 
	\end{proof}

	\section{Categories of support data}
	
	\begin{definition}\label{def:catofsupp}
		The \textit{category $\catS(\catK)$ of support data on $\catK$} is the following data:
		\begin{itemize}
			\item The objects of $\catS(\catK)$ are support data $(Y, \sigma)$;
			\item The morphisms of $\catS(\catK)$ are morphisms of support data $f\colon (Y, \sigma) \to (Z, \tau)$.
		\end{itemize}
		The category $\catS(\catK)$ has a full subcategory $\catS^\otimes(\catK)$, the \textit{category of multiplicative support data on $\catK$}. 
	\end{definition}
	
	By \cite[Theorem 4.2.2]{NVY22}, the final object of $\catS(\catK)$ is \textit{the Balmer spectrum} of $\catK$, $\Spc(\catK)$. Allowing $(\emptyset, \id)$ to be a valid support datum endows $\catS(\catK)$ with an initial object. We use the notation $\catS(-)$ to suggest functoriality. For $\catS^\otimes(-)$, this is indeed the case. 
	
	\begin{theorem}\label{thm:functorial}
		We have a contravariant functor $\catS^\otimes(-)\colon  \mathbf{mon}_{\Delta}\op \to \mathbf{Cat}$ from the category of essentially small monoidal triangulated categories to the category of essentially small categories. It is defined as follows:
		\begin{itemize}
			\item An essentially small monoidal triangulated category $\catK$ is mapped to its category of multiplicative support data $\catS^\otimes(\catK)$;
			\item An exact $\otimes$-functor $F\colon \catK \to \catL$ is mapped to the following covariant pullback functor $F^*\colon \catS^\otimes(\catL) \to \catS^\otimes(\catK)$:
			\begin{itemize}
				\item A multiplicative support theory $(Y, \sigma)$ for $\catL$ is sent to $(Y, \sigma\circ F)$.
				\item A morphism of multiplicative support theories $f\colon (Y, \sigma) \to (Z, \tau)$ for $\catL$ is mapped to the assignment \[F^*(f)\colon  (Y, \sigma\circ F) \to (Z, \tau \circ F),\] which as a map of topological spaces, sends $f$ to itself. Explicitly, we have for all objects $x \in \catK$, we have \[F^*(f)\inv(\sigma(F(x))) = \tau(F(x)).\]
			\end{itemize}
		\end{itemize}
	\end{theorem}
	\begin{proof}
		First we verify $F^*$ is a functor. We must check $(Y, \sigma\circ F)$ is a multiplicative support theory for $\catK$; this is routine and follows from $F$ being an exact $\otimes$-functor. We sketch out the tensor product property: \[ (\sigma \circ F)(x \otimes y) = \sigma(F(x)\otimes F(y)) = \sigma(F(x)) \cap \sigma(F(y)) = (\sigma \circ F)(x) \cap (\sigma \circ F)(y).\] We must verify that $F^*(f)$ is a morphism of support data, that is, for all $x \in \catK$, $F^*(f)\inv(\tau \circ F)(x) = (\sigma\circ F)(x)$, however this follows easily via the construction. Finally, it is clear that $F^*(\id) = \id$ and for $f\colon(Y, \sigma) \to (Z, \tau)$ and $g\colon(Z, \tau) \to (Z',\tau')$ morphisms of support data, that $F^*(g) \circ F^*(f) = F^*(g \circ f)$, hence $F^*$ is functorial, and thus the functor $\catS^\otimes(-)$ is well-defined. Functoriality of $\catS^\otimes(-)$ amounts to showing $\id^* = \id$ and $(G\circ F)^* = F^* \circ G^*$ which is straightforward. 
	\end{proof}
	
	On the other hand, $\catS(-)$ is not necessarily functorial. We give an explicit counterexample in Theorem \ref{ex:nocps}. However, $\catS(-)$ is functorial when restricted to essentially surjective functors. 
	
	\begin{theorem}
		Let $\catM$ denote the subcategory of $\mathbf{mon}_{\Delta}$ containing all objects, only essentially surjective morphisms, and all 2-morphisms between essentially surjective morphisms. We have a contravariant functor $\catS(-)\colon\catM\op \to \mathbf{Cat}$ as follows:
		\begin{itemize}
			\item An essentially small monoidal triangulated category $\catK$ is mapped to its category of support data $\catS(\catK)$;
			\item An essentially surjective exact $\otimes$-functor $F\colon\catK \to \catL$ is mapped to the following pullback functor $F^*\colon\catS^\otimes(\catL) \to \catS^\otimes(\catK)$:
			\begin{itemize}
				\item A multiplicative support theory $(Y, \sigma)$ for $\catL$ is sent to $(Y, \sigma\circ F)$.
				\item A morphism of multiplicative support theories $f\colon(Y, \sigma) \to (Z, \tau)$ for $\catL$ is mapped to the assignment \[F^*(f)\colon(Y, \sigma\circ F) \to (Z, \tau \circ F),\] which as a map of topological spaces, sends $f$ to itself. Explicitly, we have for all objects $x \in \catK$, we have \[F^*(f)\inv(\sigma(F(x))) = \tau(F(x)).\]
			\end{itemize}
		\end{itemize}
	\end{theorem}
	\begin{proof}
		The proof follows essentially the same as the proof of Theorem \ref{thm:functorial}, with the only major difference being the verification of axiom (e) of Definition \ref{def:supportdata}. We have:
		\begin{align*}
			\bigcup_{z \in \catK}(\sigma \circ F)(x \otimes z \otimes y) &= \bigcup_{z \in \catK}\sigma(F(x) \otimes F(z) \otimes F(y))\\
			&= \bigcup_{z' \in \catL}\sigma(F(x) \otimes z' \otimes F(y))\\
			&=\sigma(F(x)) \cap \sigma(F(y))\\
			&= (\sigma\circ F)(x) \cap (\sigma\circ F)(y),
		\end{align*}
		where equality holds between the first and second line due to essential surjectivity of $F$. 
	\end{proof}
	
	The key obstruction to complete functoriality of $\Spc$ is condition (e), as one can see in the above verification and Theorem \ref{ex:nocps}. Before moving on, we provide some commentary on our choices of support axioms and alternative approaches towards support in noncommutative tensor-triangular geometry. 
	
	\begin{remark}
		In the context of this paper, we take multiplicative support to mean support data for which the tensor product property holds, but different avenues have been explored. For instance, Negron and Pevtsova in \cite{NP23} consider ``multiplicative'' support data in the context of support for finite tensor categories under a weaker condition: the tensor product property holds only for certain objects that \textit{centralize the simple modules}. Interestingly, in most of these cases, these support theories end up having the tensor product property in practice, suggesting that perhaps verifying multiplicity can be reduced to checking multiplicativity on a small class of objects instead. 
		
		Meanwhile, the abstract question of categories whose (noncommutative) Balmer spectra have the tensor product property is the main focus of \cite{NVY22b}, and the question of when various cohomological support varieties for Hopf algebras have the tensor product property has been a subject of investigation in \cite{PW18}, \cite{BPW25}. In \cite{BPW24} and \cite{HV25} the \textit{crossed product categories} of monoidal triangulated categories with a group action were exhibited to not have the tensor product property. These categories can be thought of as semidirect products via the group action; we discuss these categories in greater detail and evaluate their complete prime spectra in Section \ref{sec:crossedprod}.
		
	\end{remark}
	
	\begin{remark}\label{rem:long}
		Similar to how in tensor-triangular geometry, the standard support axioms break down for ``big'' tensor-triangulated categories (as addressed by Barthel, Heard, and Sanders \cite{BHS23} building on Balmer-Favi support), the question of what support axioms are ``good'' for a small monoidal triangulated category is a more open-ended question. To understand different notions of support arising from different choices of support axioms, we suggest the following approach. Given a list of support axioms for $\catK$, we define $\catS(\catK)$ in the same manner as Definition \ref{def:catofsupp}, that is, all pairs $(Y, \sigma)$ of support data on $\catK$ that satisfy the support axioms. Given a list of support axioms, the fundamental questions one may ask are: 
		
		\begin{itemize}
			\item Is $\catS(-)$ functorial?
			\item Does $\catS(\catK)$ have a final object, i.e. a spectrum $\Spc(\catK)$?
		\end{itemize}
		
		Alternatively, given an abstract support theory $\Spc$ for every essentially small monoidal triangulated category $\catK$, the fundamental questions are:
		
		\begin{itemize}
			\item Is $\Spc(\catK)$ functorial? 
			\item Is $\Spc(\catK)$ the final object in a category $\catS(\catK)$ for some list of support axioms?
		\end{itemize}
		
		Once these are established, the deeper questions are:
		
		\begin{itemize}
			\item What do the specialization closed/Thomason subsets of $\Spc(\catK)$ parametrize? Can this be read off purely from the support axioms, even without having an explicit description of the Balmer spectrum? 
			\item Is $\Spc(\catK)$ large enough to be informative, but small enough to be computable? 
			\item Given an exact $\otimes$-functor $F\colon  \catK \to \catL$, when is the induced map on Balmer spectra injective or surjective? 
		\end{itemize}
		
		One must wonder if there is a all-encompassing notion of support. In the ring-theoretic case, Reyes \cite{Rey12} determines that, if one wants the spectrum to be functorial and coincide with commutative ring spectra for commutative rings, one cannot guarantee the spectrum will always be nonempty. Similarly, we want that our support axioms should always coincide with the usual (multiplicative) support axioms for symmetric or braided monoidal categories. We ask the following question: does a universal, functorial, necessarily non-empty support theory exist? The final section lays out the groundwork to solve this answer in the same manner as Reyes, but at present the answer remains unclear. 
	\end{remark}
	
	\begin{remark}
		We note one downside of working with the complete prime spectrum or more general support theories: much of the analysis of the Balmer spectrum (as a set of primes) and in particular, the specialization-closed subsets, utilizes frame-theoretic considerations pioneered in \cite{BKS07} and \cite{KP16}. See \cite{BKS20}, \cite[Appendix A]{NVY23}, \cite{MR23}, \cite{Kra24} for further examples. These techniques work when the poset $\Spc(\catK)$ is a \textit{(semi)lattice}, for instance when $\Spc(\catK)$ is the Balmer spectrum as specified in \cite{NVY22}, the set of two-sided thick $\otimes$-ideals. In this setting, one usually has an abstract ``prime lifting'' theorem - such a theorem may not hold for the complete prime spectrum of a monoidal triangulated category. Such considerations are amenable for classification questions, see for instance \cite{K23}.
	\end{remark}
	
	We have presented here a general approach to evaluating different axioms of support, but for the context of this paper we reserve $\catS(\catK)$ and $\catS^\otimes(\catK)$ for the support axioms defined in Definition \ref{def:catofsupp}, and mainly focus on understanding multiplicative or functorial support data. Throughout this article, we will aim to answer the above questions.
	
	\section{Partial functoriality of the noncommutative Balmer spectrum}

	\begin{remark}
		Why is functoriality important, besides simply being a nice property to have? We turn to the classical symmetric monoidal category case for inspiration. Often when the (commutative) Balmer spectrum cannot be discerned by standard means such as comparison maps as in \cite{Bal10} or \cite{NVY24}, it is discerned from ``fiber functors,'' a family of functors which together detect vanishing or tensor-nilpotence of morphisms (certain techniques allow for weaker assumptions than conservativity, see \cite[Definition 1.1]{BCHS24}). For instance, Balmer and Gallauer in \cite{BG23} construct modular fixed points functors to recover $\Spc(\on{K}_b(p\on{-perm}(kG)))$ as a set, where $\on{K}_b(p\on{-perm}(kG))$ denotes the bounded homotopy category of $p$-permutation $kG$-modules. Similarly, in \cite{BS17}, Balmer and Sanders use geometric fixed points functors to deduce $\Spc(\on{SH}(G)^c)$, where $\on{SH}(G)^c$ denotes the compact part of the $G$-equivariant stable homotopy category. Conservativity techniques have been developed in detail for symmetric monoidal categories in \cite{Bal18}, \cite{BCHS24} and recently \cite{BG25}.
	\end{remark}
	
	We confirm the expectation that $\Spc(\catK)$ is generally non-functorial in Theorem \ref{ex:nocps}. Of course, this runs counter to the well-known property that $\Spc$ is functorial for symmetric tensor-triangulated categories, \cite[Proposition 3.6]{Bal05}. However in some special circumstances, $\Spc$ induces a continuous homomorphism of topological spaces, such as in \cite[Proposition 2.1.3]{Vas24}. Similarly to $\catS(-)$ case, when we restrict to essentially surjective functors, we obtain functoriality. 
	
	\begin{prop}\label{prop:surjmaponspectra}
		Let $F\colon\catK \to \catL$ be a exact $\otimes$-functor. If $F$ is essentially surjective, then $\varphi := \Spc(F)\colon\Spc(\catL) \to \Spc(\catK),\catQ \mapsto F\inv(\catQ)$ is a continuous, well-defined homomorphism satisfying $\varphi\inv(\supp_\catK(x)) = \supp_\catL(F(x))$ for all $x \in \catK$. 
	\end{prop}
	\begin{proof}
		Suppose $\catQ \in \Spc(\catL)$ is a prime thick $\otimes$-ideal. Verification that $F\inv(\catQ)$ is a thick $\otimes$-ideal is straightforward, we show it is prime. Suppose for fixed $x,y\in\catK$ that $x \otimes k \otimes y \in F\inv(\catQ)$ for all $k \in \catK$. Then $F(x) \otimes F(k) \otimes F(y) \in \catQ$ for all $k \in \catK$, and since $F$ is essentially surjective,  equivalently $F(x) \otimes l \otimes F(y) \in \catQ$ for all $l \in \catL$. Since $\catQ$ is prime, either $F(x)$ or $F(y)$ belongs to $\catQ$, hence $x$ or $y$ belongs to $F\inv(\catQ)$. Since $F$ is a monoidal functor, $\bbone\not\in F\inv(\catQ)$, thus $\varphi$ is well-defined. Since $\Spc(\catK)$ has a basis of open sets $U(x)$ for all $x \in \catK$, it suffices to verify the support formula for $\varphi$ holds. We verify:
		\begin{align*}
			\supp_\catL(F(x)) &= \{\catQ \in \Spc(\catL) \mid F(x) \not\in \catQ\}\\
			&= \varphi\inv(\{\catP \in \Spc(\catK)\mid x \not\in \catP\})\\
			&= \varphi\inv(\supp_\catK(x)).
		\end{align*}
	\end{proof}
	
	Given any thick $\otimes$-ideal $\catI$ of $\catK$, one can perform the localization at $\catI$, i.e. taking the Verdier quotient $\catK/\catI$, which is again monoidal triangulated. We have a canonical localization functor $\catK \twoheadrightarrow \catK/\catI$, which is surjective on objects. 
	
	\begin{corollary}
		Let $\catI$ be a thick $\otimes$-ideal of $\catK$, and let $q\colon \catK \to \catK/\catI$ be a localization. The map $\varphi:=\Spc(q)\colon \Spc(\catK/\catI) \to \Spc(\catK)$ induces a homeomorphism between $\Spc(\catK/\catI)$ and the subspace $\{\catP \in \Spc(\catK) \mid \catI \subseteq \catP\}$ of $\Spc(\catK)$ of prime thick $\otimes$-ideals containing $\catI$. Similarly, the map $\varphi\cp: \Spc\cp(q)\colon \Spc\cp(\catK/\catI) \to \Spc(\catK)$ induces a homeomorphism between $\Spc(\catK/\catI)$ and the subspace $\{\catP \in \Spc\cp(\catK) \mid \catI \subseteq \catP\}$ of $\Spc\cp(\catK)$. 
	\end{corollary}
	\begin{proof}
		We show the statement for the noncommutative Balmer spectrum; the latter statement follows analogously. Well-definedness of $\varphi := \Spc(q)$ is guaranteed by the previous proposition. Set $V := \{\catP \in \Spc(\catK) \mid \catI \subseteq \catP\}$. For any $\catQ \in \Spc(\catK/\catI)$, we have $\varphi(\catQ)= q\inv(\catQ) \supset q\inv(0) = \catI$, or in other words, $\im(\varphi) \subseteq V$. Proposition \ref{prop:injmaponspec} asserts $\varphi$ is injective, since $q$ is surjective. Conversely. if $\catP \in \Spc(\catK)$ contains $\catI$, then it is straightforward to see $q(\catP)$ is a prime thick $\otimes$-ideal, and $q\inv(q(\catP)) = \catP$, so $\varphi$ induces a continuous bijection $\Spc(\catK/\catI)\to V$. 
		
		Now, let $x \in \catK$ and $\catP \in V$. We have $x \in \catP$ if and only if $x \in q(\catP)$ (with $x$ regarded as an object of $\catK/\catI$), so $\varphi(\supp_{\catK/\catI}(x)) = \supp_{\catK}(x) \cap V$, hence $\varphi: \Spc(\catK/\catI) \to V$ is a closed map, thus a homeomorphism as desired. 
	\end{proof}
	
	We can in fact give a criterion for when an exact $\otimes$-functor $F\colon \catK \to \catL$ induces a well-defined map on noncommutative Balmer spectra. 
	
	\begin{theorem}\label{cor:wheniswelldef}
		Let $F\colon \catK \to\catL$ be an exact $\otimes$-functor, and denote its image by $F(\catK)\subseteq \catL$. Then $F$ induces a continuous map on Balmer spectra $\varphi := \Spc(F): \Spc(\catL) \to \Spc(\catK)$ given by $\catQ \mapsto F\inv(\catQ)$ if and only if for all prime thick $\otimes$-ideals $\catQ \in \Spc(\catL)$, the intersection $F(\catK) \cap \catQ$ satisfies the primality property, \[\text{For all $a, b \in F(\catK)$, } a \otimes F(\catK) \otimes b \subseteq F(\catK) \cap \catQ \text{ implies } a \in \catQ \text{ or } b\in \catQ.\] 
	\end{theorem}
	\begin{proof}
		It is routine to verify that for every thick $\otimes$-ideal $\catJ \subseteq \catL$, $F\inv(\catJ)$ is a thick $\otimes$-ideal of $\catK$. Therefore, for $\Spc(F)$ is well-defined if and only if the primality property is preserved for all primes $\catQ \in \Spc(\catL)$. We have that $F$, as a monoidal (but not necessarily triangulated, since $F(\catK)$ may not be triangulated if $F$ is not full) functor factors as the surjective monoidal functor $F\colon  \catK \to F(\catK)$ composed with the inclusion $\iota: F(\catK) \subseteq \catL$. The same argument as Proposition \ref{prop:surjmaponspectra} demonstrates that for every $\catQ \in \Spc(\catL)$, $F\inv(F(\catK) \cap \catQ)\subseteq \catK$ satisfies the tensor product property by surjectivity of $F\colon \catK \to F(\catK).$ We turn to the inclusion $\iota$.
		
		If $\iota\inv$ preserves the primality property for all intersections of primes $\catQ \cap F(\catK)$ with $\catQ \in \Spc(\catL)$, it follows that the composition $(\iota\circ F)\inv = F\inv\circ \iota\inv$ preserves the property as well, so $\Spc(F)$ is well-defined. Conversely if the property is not satisfied, then there exist $a,b \in F(\catK)$ such that $a \otimes F(\catK) \otimes b \in F(\catK) \cap \catQ$ but $a, b \not\in \catQ$. Therefore, there exist $a', b' \in \catK$ with $F(a') = a$, $F(b') = b$ such that $a' \otimes \catK \otimes b' \subseteq F\inv(\catQ) = F\inv(\catQ \cap F(\catK))$, but $a', b'\not\in F\inv(\catQ)$, hence $\Spc(F)$ is not well-defined. Continuity follows the same as in Proposition \ref{prop:surjmaponspectra}.
	\end{proof}
	
	We obtain a standard tensor-triangular geometric result for noncommutative tensor-triangular geometry: idempotent completion does not affect the Balmer spectrum. This is a noncommutative version of \cite[Proposition 3.13]{Bal05}.
	
	\begin{theorem}
		Let $\catL$ be a monoidal triangulated category and let $\catK \subseteq \catL$ be a full monoidal triangulated subcategory with the same unit and which is cofinal, i.e. for any object $x \in \catL$, there exists $y \in \catL$ such that $x \oplus y \in \catK$. Then the map $\catQ \mapsto \catQ\cap \catK$ defines a homeomorphism $\Spc(\catL) \xrightarrow{\cong} \Spc(\catK)$. Moreover, this restricts to a homeomorphism $\Spc\cp(\catL) \xrightarrow{\cong} \Spc\cp(\catK)$.
	\end{theorem}
	\begin{proof}
		First, note we may replace $\catK$ with its isomorphic-closure, so we can assume $\catK$ is replete. The map $\Spc(\catL) \to \Spc(\catK)$, if it is well-defined, is simply $\Spc(\iota)$, where $\iota: \catK \hookrightarrow \catL$ denotes inclusion. To see $\Spc(\iota)$ is well-defined, for any $\catQ \in \Spc(\catL)$, $x \otimes (z \oplus z') \otimes y \in \catQ$ implies $x \otimes z \otimes y \in \catQ$ and $x \otimes z' \otimes y \in \catQ$ since $\catQ$ is thick. Therefore by cofinality, for any $x, y \in \catK$, $x \otimes \catK \otimes y \subseteq \catQ$ implies $x \otimes \catL \otimes y \subseteq \catQ$ implies $x \in \catQ$ or $y \in \catQ$, so $\catQ \cap \catK$ is indeed a prime thick $\otimes$-ideal of $\catK$, and $\Spc(\iota)$ is indeed a well-defined map on Balmer spectra. 
		
		Now, recall that for any $x \in \catL$, we have $x \oplus \Sigma x \in \catK$ (c.f. the proof of \cite[Proposition 3.13]{Bal05}). Moreover, if $x \oplus x'$ belongs to a triangulated subcategory of $\catK$ (e.g. a prime), then so does $x \oplus \Sigma x$. Therefore, given a prime $\catP \in \Spc(\catK)$, we have an equality of full subcategories \[\widetilde{\catP}:= \{x \in \catL \mid x \oplus \Sigma x \in \catP\} = \{x \in \catL \mid \exists x' \in \catL \text{ such that } x \oplus x' \in \catP\}.\] We claim $\widetilde{\catP}$ is a prime of $\catL$. It is routine to verify $\widetilde{\catP}$ is a thick $\otimes$-ideal. To verify it is prime, suppose $x \otimes z \otimes y \in \widetilde{\catP}$ for all $z \in \catL$, and suppose $x \not\in \widetilde{\catP}.$ Set $w := x \oplus \Sigma x$, then by construction, $w \in \catK \setminus \catP$, but \[w \otimes z \otimes y \cong (x \otimes z \otimes y) \oplus \Sigma (x \otimes z \otimes y) \in \catP\] for all $z \in \catK$, therefore \[w \otimes z \otimes (y \oplus \Sigma y) \cong (w \otimes z \otimes y) \oplus\Sigma(w\otimes z\otimes y) \in \catP\]
		for all $z \in \catK$, hence $y \oplus \Sigma y \in \catP$ by primality of $\catP$. Thus, $y \in \widetilde{\catP}$, as desired. 
		
		It follows by $\catP$ thick that $\widetilde{\catP} \cap \catK = \catP$, so the assignment $\catP \mapsto \widetilde{\catP}$ is a right inverse to $\Spc(\iota)$. On the other hand, if $\catQ \in \Spc(\catL)$, we claim $\catQ = \widetilde{\catP}$ where $\catP := \catQ \cap \catK$. Indeed, if $x \in \catQ$, then $x \oplus \Sigma x \in \catK$, hence $x \oplus \Sigma x \in \catQ$, so $\catQ \subseteq \widetilde{\catP}$. Conversely, $\widetilde{\catP} \subseteq \catQ$ by thickness of $\catQ$. So $\Spc(\iota)$ is a continuous bijection with inverse given by $\catP \mapsto \widetilde{\catP}$. 
		
		Now, for any $x \in \catL$, and $\catP\in \Spc(\catK)$, we have $x \in \widetilde{\catP}$ if and only if $x \oplus \Sigma x \in \widetilde{\catP}$ if and only if $x \oplus \Sigma x \in \catP$. Therefore, the image of $\Spc(\iota)$ of $\supp_\catL(x)$ is $\supp_\catK(x \oplus \Sigma x),$ so $\Spc(\iota)$ is a closed map, and thus a homeomorphism, as desired. 
		
		For the final statement, one needs to verify that $\Spc(\iota)$ and the assignment $\catP \mapsto \widetilde{\catP}$ for $\catP \in \Spc(\catK)$ preserve complete primality. The former follows from functoriality of $\Spc\cp$, proven in Theorem \ref{thm:thm5}. For the latter, suppose $x \otimes y \in\widetilde{\catP}$ and $x \not\in \widetilde{\catP}$. Let $z := x \oplus \Sigma x$, then we have $z \in \catK \setminus \catP$ and $z \otimes y \cong (x \otimes y) \oplus \Sigma (x \otimes y) \in \catP$. Therefore, $z \otimes (y \oplus \Sigma y) \cong (z \otimes y) \oplus \Sigma(z \otimes y) \in \catP$, so $y \oplus \Sigma y \in \catP$, hence $y \in \widetilde{\catP}$. 
	\end{proof}
	
	\begin{corollary}
		Let $\catK^\natural$ denote the idempotent completion of $\catK$, and let $\iota: \catK \hookrightarrow \catK^\natural$ denote the canonical inclusion. The induced map on the noncommutative Balmer spectrum $\Spc(\iota): \Spc(\catK^\natural) \to \Spc(\catK)$ is a well-defined homeomorphism, which restricts to a homeomorphism $\Spc\cp(\iota): \Spc\cp(\catK^\natural) \to \Spc\cp(\catK)$. 
	\end{corollary}
	
	\section{Functoriality and universality of the complete prime spectrum}
	
	In this section, we verify the key properties of the complete prime spectrum: universality and functoriality. 
	
	\begin{lemma}{\cite[Lemma 3.3]{Bal05}}\label{lem:equalityoffuncs}
		Let $S$ be a set, assume $\Spc^{cp}(\catK) \neq \emptyset$, and let $f_1, f_2\colon S \to \Spc^{cp}(\catK)$ be two maps such that $f_1\inv(\supp^{cp}(x)) = f_2\inv(\supp^{cp}(x))$ for any $x \in \catK$. Then $f_1 = f_2$. 
	\end{lemma}
	\begin{proof}
		Observe that for any $x \in \catK$ and $s \in S$, by assumption we have $f_1(s) \in \supp^{cp}(x)$ if and only if $f_2(s) \in \supp^{cp}(x)$. Fix $s \in S$, then the following two closed subsets of $\Spc(\catK)$ coincide: \[\bigcap_{f_1(s) \in \supp^{cp}(x)} \supp^{cp}(x) = \bigcap_{f_2(s) \in \supp^{cp}(x)} \supp^{cp}(x).\] However the left hand side is $\overline{\{f_1(x)\}}$ and the right hand side is $\overline{\{f_2(x)\}}$, hence Proposition \ref{prop:t0} implies $f_1(x) = f_2(x)$ for all $s \in S$, so $f_1 = f_2$ as desired.  
	\end{proof}
	
	Now we prove, as promised, that $\Spc^{cp}(\catK)$ is universal: it is the final object of $\catS^\otimes(\catK)$. 
	
	\begin{theorem}\label{thm:thm6}
		The complete prime spectrum $(\Spc^{cp}(\catK), \supp^{cp})$ is the final multiplicative support data on $\catK$, i.e. the final object of $\catS^\otimes(\catK)$. That is, if it is nonempty, $(\Spc^{cp}(\catK), \supp^{cp})$ is a multiplicative support data, and for any other multiplicative support data $(Y, \sigma)$ on $\catK$, there exists a unique continuous map $f\colon Y \to \Spc^{cp}(\catK)$ such that $\sigma(x) = f\inv(\supp^{cp}(x))$. Explicitly, the map is defined by \[f(x) = \{a \in \catK \mid x \not\in \sigma(a)\}.\]
		On the other hand, if $\Spc\cp(\catK) = \emptyset$, then $(\Spc\cp(\catK), \supp)$ is the only multiplicative support datum on $\catK$, therefore it is trivially the final object of $\catS^\otimes(\catK)$. 
	\end{theorem}
	
	\begin{proof}
		First we consider when $\Spc\cp(\catK)$ is nonempty. We have previously verified that $(\Spc^{cp}(\catK), \supp^{cp})$ is a multiplicative support data. To verify the universal property, suppose $(Y,\sigma)$ is a multiplicative support data on $\catK$. Uniqueness of $f$ follows from Lemma \ref{lem:equalityoffuncs}. To check that the map $f$ as defined is the unique map, it is a routine verification that $f(x)$ is a thick $\otimes$-ideal. To show $f(x)$ is completely prime, let $a \otimes b \in f(x)$. This means that $x \not\in \sigma(a\otimes b) = \sigma(a) \cap \sigma(b)$, therefore $x \not\in \sigma(a)$ or $x\not\in\sigma(b)$, i.e. $a \in f(x)$ or $b \in f(x)$, as desired. Now, by definition, we have that $f(x) \in \supp^{cp}(a) $ if and only if $a \not\in f(x)$ if and only if $x \in \sigma(a)$. Therefore, $f\inv(\supp(a)) = \sigma(a)$, so $f$ is indeed the unique continuous map specified. Finally, this map is continuous by construction, since the compliments of $\supp^{cp}(a)$ give a basis of open sets of $\Spc\cp(\catK)$. 
		
		Now, suppose $\Spc\cp(\catK)$ is empty, and suppose that $(Y, \sigma)$ is a multiplicative support datum. In particular, $(Y, \sigma)$ is a support datum, therefore we have a unique map $f\colon (Y, \sigma) \to (\Spc(\catK), \supp)$. Moreover, we have that $\sigma(x) = f\inv(\supp(x))$ for all $x \in \catK$. Since $\sigma$ satisfies the tensor product property, we have for all $x,y \in \catK$ the following equalities
		\begin{align*}
			f\inv(\supp(x \otimes y)) &= \sigma(x \otimes y)\\
			&= \sigma(x) \cap \sigma(y)\\
			&= f\inv(\supp(x) \cap \supp(y))\\
			&= f\inv(\supp(x)) \cap f\inv(\supp(y)).
		\end{align*}
		However, since there are no completely prime ideals, for any $\catP \in \Spc(\catK)$, there exist $x, y \in \catK$ for which $x \otimes y \in \catP$ and $x,y \not\in\catP$. Equivalently, we have $\catP \in \supp(x),\supp(y)$ but $\catP \not\in \supp(x\otimes y)$. The above chain of equalities implies $\catP \not\in \im(f)$, but this holds for all primes $\catP \in \Spc(\catK)$. Hence $\im(f) = \emptyset$, so necessarily $Y = \emptyset$, as desired. 
	\end{proof}
	
	This universality provides an alternative short ``universal'' proof of the fact that for a monoidal triangulated category, every prime is completely prime if and only if the tensor product property holds, \cite[Theorem 3.1.1]{NVY22b}.
	
	\begin{corollary}
		For every monoidal triangulated category $\catK$, the following are equivalent.
		\begin{enumerate}
			\item The Balmer spectrum $(\Spc(\catK),\supp)$ has the tensor product property.
			\item Every prime ideal of $\catK$ is completely prime.
		\end{enumerate}
	\end{corollary}
	\begin{proof}
		If every prime ideal of $\catK$ is completely prime, then $\Spc^{cp}(\catK) \hookrightarrow \Spc(\catK)$ is the identity, hence $(\Spc(\catK), \supp)$ is multiplicative. Conversely, if $(\Spc(\catK), \supp)$ has the tensor product property, then the unique map map $f:\Spc(\catK) \to \Spc^{cp}(\catK)$ defined by $f(\catP) = \{a \in \catK \mid \catP \not\in \supp(a)\}$ is clearly the identity. Thus $(\Spc(\catK),\supp) = (\Spc^{cp}(\catK),\supp^{cp})$, so $(\Spc(\catK),\supp)$ has the tensor product property. 
	\end{proof}

	\begin{corollary}
		Let $\catK$ be a monoidal triangulated category satisfying the tensor product property. Then every support data and weak support data (\cite[Definition 4.4.1]{NVY22}) $(Y, \sigma)$ satisfies the tensor product property.
	\end{corollary}
	\begin{proof}
		We have that $\Spc(\catK) = \Spc\cp(\catK)$. Now, let $(Y, \sigma)$ be a support data or weak support data for $\catK$. By universality, we have a unique map $f: (Y, \sigma) \to (\Spc(\catK) = \Spc\cp(\catK), \supp)$ satisfying $\sigma(x) = f\inv(\supp(x))$ for all $x \in \catK$. Now, we have \[\sigma(a \otimes b) = f\inv(\supp(a \otimes b)) = f\inv(\supp(a) \cap \supp(b)).\] On the other hand, \[\sigma(a) \cap \sigma(b) = f\inv(\supp(a)) \cap f\inv(\supp(b)) = f\inv(\supp(a) \cap \supp(b)),\] hence $\sigma(a) \cap \sigma(b) = \sigma(a \otimes b)$, as desired. 
	\end{proof}
	
	Finally, we show $\Spc^{cp}$ is fully functorial, as promised. 
	
	\begin{theorem}\label{thm:thm5}
		Given a monoidal $\otimes$-functor $f\colon \catK \to \catL$, the map \[\varphi:= \Spc^{cp}(F)\colon  \Spc^{cp}(\catL) \to \Spc^{cp}(\catK),\quad \catQ \mapsto F\inv(\catQ)\] is well-defined, continuous, and for all objects $x \in \catK$, we have \[\varphi\inv\big(\supp^{cp}_{\catK}(x)\big) = \supp^{cp}_\catL(F(x)).\] This defines a contravariant functor $\Spc^{cp}$ from the category of essentially small monoidal triangulated categories to the category of topological spaces. 
	\end{theorem}
	
	\begin{proof}
		First we show well-definedness. It is easy to check that $F\inv(\catP)$ is a thick $\otimes$-ideal. For primality, suppose $a \otimes b \in F\inv(\catP)$, then $F(a) \otimes F(b) \in \catP$. By primality, either $F(a)$ or $F(b)$ belongs to $\catP$, hence $a$ or $b$ belongs to $F\inv(\catP)$ by definition. Note that this implies that if $\Spc\cp(\catL)$ is nonempty, then $\Spc\cp(\catK)$ necessarily is as well. On the other hand, if $\Spc\cp(\catL)$ is empty, this poses no problems since the empty set is the initial object in the category of topological spaces. 
		
		We verify the equality $(\Spc F)\inv\big(\supp^{cp}_{\catK}(x)\big) = \supp^{cp}_\catL(F(x))$. By definition, \[\supp^{cp}_\catL(F(x)) = \{\catQ \in \Spc\cp(\catL) \mid F(x) \neq \catQ\}.\] On the other hand, \[(\Spc F)\inv\big(\supp^{cp}_{\catK}(x)\big) = \{\catQ \in \Spc\cp(\catL) \mid  x\not\in (\Spc(F))(\catQ)\}.\] By definition, $x\not\in (\Spc(F))(\catQ)$ if and only if $x \not\in F\inv(\catQ)$ if and only if $F(x) \not\in \catQ$, so the equality holds. 
		
		Since the compliments of supports form a basis of the topology of the complete prime spectrum, it follows that $\Spc(F)$ is continuous. Finally, $\Spc^{cp}(G \circ F) = \Spc^{cp}(F) \circ \Spc^{cp}(G)$ is straightforward. 
	\end{proof}
	
	We make an easy observation: functoriality provides an easy check for when $\Spc\cp(\catK) \neq \emptyset$, the existence of a ``fiber functor'' to a more well-understood category with nonempty complete prime spectrum. 
	
	\begin{corollary}
		Let $F\colon  \catK\to \catL$ be a (necessarily nonzero) exact $\otimes$-functor and suppose $\Spc\cp(\catL) \neq \emptyset$. Then $\Spc\cp(\catK) \neq \emptyset$.
	\end{corollary}
	\begin{proof}
		We prove the contrapositive. If $\Spc\cp(\catK) =\emptyset$ then we have a morphism of topological spaces $\Spc\cp(F): \Spc\cp(\catL) \to \Spc\cp(\catK) = \emptyset$. The only morphism into $\emptyset$ is the identity $\id_\emptyset: \emptyset\to\emptyset$, thus $\Spc\cp(\catL) = \emptyset$. 
	\end{proof}

	We also have an easy extension of \cite[Theorem A]{Vas24} - another instance when $\Spc$ is functorial. First, we recall the Drinfeld centralizer and center of a tensor category, following \cite{Vas24}. For $\catD \subseteq \catC$ a tensor subcategory, we can associate the \textit{Drinfeld centralizer} $Z^{\catD}(\catC)$ of $\catD$ in $\catC$, the finite tensor category of pairs $(V,\gamma_V)$, where $V\in \catC$ is an object of $\catC$ and $\gamma_V\colon V\otimes - \to - \otimes V$ is a natural isomorphism of functors satisfying the constraint \[(X \otimes \gamma_{V,Y}) \circ (\gamma_{V,X} \otimes Y )= \gamma_{V, X \otimes Y}\] for all $X, Y \in \catD$. Then $Z^\catD(\catC)$ is the category of pairs $(V, \gamma_V)$ for an object $V \in \catC$ and a choice of $\catD$-centralizing structure $\gamma_V$ on $V$. Note this choice is in general non-canonical, and thus even is $\catC$ is triangulated and $\catD$ is a triangulated subcategory, $Z^\catD(\catC)$ may not be triangulated. However, we obtain a forgetful functor $Z^\catD(\catC) \to \catC$, clearly injective on objects. 
	
	In the case of $\catD = \catC$, we simply write $Z(\catC)$ and obtain the \textit{Drinfeld center} of $\catC$. The other case is when $\catD$ is the tensor subcategory generated by the semisimple objects of $\catC$.
	
	\begin{prop}\label{prop:injmaponspec}
		Let $\catC,\catD$ be a finite tensor categories, and let $F\colon \stab(\catC)\to\stab(\catD)$ be a monoidal triangulated functor. If the essential image of $F$ is contained in the image of the forgetful functor $\stab(Z(\catD)) \to \stab(\catD)$, where $Z(\catD)$ denotes the Drinfeld center then $\Spc(F)\colon  \Spc(\stab(\catD)) \to \Spc(\stab(\catC))$ is well-defined.
	\end{prop}
	\begin{proof}
		In this case, $F$ factors as $f\colon \stab(\catC) \to \stab(Z(\catD)) \hookrightarrow \stab(\catD)$. The latter functor induces a continuous map on spectra by \cite[Proposition 2.1.3]{Vas24}, and $\stab(Z(\catD))$ is braided, therefore $\Spc\cp(\stab(Z(\catD))) =\Spc(\stab(Z(\catD)))$, hence we obtain a corresponding continuous homomorphism \[\Spc(\stab(Z(\catD))) = \Spc\cp(\stab(Z(\catD)))\to \Spc\cp(\stab(\catC)) \hookrightarrow \Spc(\stab(\catC)).\]
	\end{proof}
	
	As $\Spc\cp$ is functorial, we can ask when is the induced topological map $\varphi = \Spc^{cp}(F)\colon  \Spc^{cp}(\catL) \to \Spc^{cp}(\catK)$ injective or surjective. One case is easy.
	
	\begin{prop}
		If $f\colon \catK \to \catL$ is an essentially surjective monoidal triangulated functor, then $\Spc(F)$ and $\Spc\cp(F)$ are injective.
	\end{prop}
	\begin{proof}
		In this case, any (complete) prime $\catQ \subset \catL$ satisfies that $\langle F(F\inv(\catQ))\rangle = \catQ$ (if $F$ is surjective on objects, one can remove the brackets), hence $F\inv(\catQ_1) = F\inv(Q_2)$ forces $\catQ_1 = \catQ_2$. 
	\end{proof}
	
	Surjectivity will be explored in Section \ref{sec:surj}.
	
	\section{Radical thick tensor-ideals}
	
	In Remark \ref{rem:long}, we claimed that any reasonable choice of support theory axioms should satisfy that the specialization-closed subsets of its Balmer spectrum should parametrize something about our category $\catK$. Indeed, perhaps unsurprisingly we will see in this section that $\Spc\cp(\catK)$ parametrizes \textit{radical} thick $\otimes$-ideals. 
	
	% We'll investigate the specialization closed subsets for multiplicative support data, i.e. the complete prime spectrum $\Spc\cp(\catK)$. Perhaps unsurprisingly, this verification mirrors the classical setting of \cite{Bal05}, where \textit{radical} thick $\otimes$-ideals come into play.  
	
	\begin{definition}
		Let $\catI$ be a thick $\otimes$-ideal of $\catK$. We say $\catI$ is \textit{radical} if $\sqrt\catI = \catI$ or equivalently, for all $x \in \catK$, $x \in \catI$ if and only if $x^{\otimes n} \in \catI$ for all $n \geq  1$.
	\end{definition}
	
	In fact, complete primes are nothing more than radical primes. 
	
	\begin{prop}\label{prop:cpprimeradical}
		Let $\catP$ be a thick $\otimes$-ideal of $\catK$. Then $\catP$ is completely prime if and only if $\catP$ is prime and radical.
	\end{prop}
	\begin{proof}
		If $\catP$ is completely prime, then $\catP$ is prime and radical. Conversely, suppose for contradiction that $\catP$ is prime and radical, but not completely prime. Then there exist $x , y \in \catK$ for which $y \otimes x \in \catP$ but neither $x \in \catP$ nor $y \in \catP$. Primality of $\catP$ implies that there exists $z \in \catK$ for which $x \otimes z \otimes y \not\in \catP$. Then $\catP$ radical implies that $(x \otimes z \otimes y)\otimes (x \otimes z \otimes y) \not\in\catP$, but this is a contradiction since $y \otimes x \in \catP$ and $\catP$ is a $\otimes$-ideal. 
	\end{proof}
	
	Now, following \cite{Kra24}, we have an analogue of \cite[Theorem 4.10]{Bal05}. However we must be careful - rather than using the language of Thomason subsets, we instead utilize the Hochster dual topology on $\Spc\cp(\catK)$. 
	
	\begin{definition}\label{def:hochster}
		We define the \textit{Hochster dual} topology of $\Spc\cp(\catK)$, denoted \newline $\Spc\cp(\catK)^\vee$, as follows. As a set, $\Spc\cp(\catK)^\vee := \Spc\cp(\catK)$, and the supports $\supp(x)$ for $x \in \catK$ generate yield a basis of \textit{open} sets for the topology of $\Spc\cp(\catK)^\vee$.
	\end{definition}
	
	\begin{theorem}\label{thm:cpparametrization}
		We have an order-preserving bijection from the set $\catR$ of radical thick $\otimes$-ideals to the set $\catS$ of \textit{open} subsets of $\Spc\cp(\catK)^\vee$ induced by the assignment \[\catI \mapsto \bigcup_{x \in \catI} \supp\cp(x).\] In particular, every radical thick $\otimes$-ideal of $\catK$ is an intersection of completely prime thick $\otimes$-ideals of $\catK$.
	\end{theorem}
	\begin{proof}
		This essentially follows from \cite[Section 6]{Kra24}; we summarize here. For two $x, y \in \catK$, write $x \sim y$ if and only if $\langle x\rangle = \langle y \rangle$, where $\langle x\rangle$ denotes the \textit{radical} thick $\otimes$-ideal generated by $x$. Then $\sim$ is an equivalence relation, and $L(\catK, \otimes) := \on{Ob}(\catK)/\sim$ is partially ordered via inclusion. \cite[Proposition 14]{Kra24} asserts that $L(\catK, \otimes)$ is a distributive lattice, and its ideal lattice $\on{Id}(L(\catK, \otimes))$ identifies with the lattice of radical thick $\otimes$-ideals of $\catK$. The bijection of the first statement now follows from \cite[Corollary 10]{Kra24} if the prime ideals of the lattice $L(\catK, \otimes)$ are precisely the complete primes of $\catK$, i.e. we have an equality of sets $\Spc\cp(\catK) = \Spc(L(\catK, \otimes))$. Indeed, this follows from the identity $\langle x \rangle \wedge \langle y \rangle = \langle x \otimes y \rangle$ of \cite[Lemma 13]{Kra24}, so $\langle x\rangle$ is a completely prime thick $\otimes$-ideal of $\catK$ if and only if it is prime as an element of $L(\catK, \otimes)$. 
		
		The latter statement follows from the inverse of the map, which sends an open set $U$ to $\{x \in L(\catK, \otimes) \mid \supp\cp(x) \subseteq U\}$, which is nothing more than the intersection of all primes of $L(\catK, \otimes)$, i.e. complete primes, not contained in $U$ (c.f. \cite[Lemma 4.8]{Bal05}).
	\end{proof}
	
	One must be careful with the usage of Thomason subsets here, since it is not clear if every open $U\cp(x)$ of $\Spc\cp(\catK)$ is quasi-compact - however, this indeed holds if every prime is complete. In fact, when $\Spc(\catK) = \Spc\cp(\catK)$, we informally note that nearly every topological property proven in \cite[Section 2]{Bal05} holds, essentially due to the crucial prime lifting lemma Proposition \ref{prop:primelift} holding for complete primes. 
	
	\section{Intermezzo: Quasi-support data}
	
	While we are on the topic of finding the universal support datum for a given set of support axioms, we address \textit{quasi support data}, as introduced in \cite[Section 7.1]{NVY22}.
	
	\begin{definition}
		A \textit{(noncommutative) quasi support datum} for a monoidal triangulated category $\catK$ is a map \[\sigma\colon \catK \to \calX_{cl}(Y)\] for a topological space $Y$ such that the following hold:
		\begin{enumerate}
			\item $\sigma(0) = \emptyset$ and $\sigma(\bbone) = Y$;
			\item $\sigma(x \oplus y) = \sigma(x) \cup \sigma(y)$;
			\item $\sigma(\sum x) = \sigma(x)$;
			\item $\sigma(a) \subseteq \sigma(b) \cup \sigma(c)$ for all distinguished triangles $a \to b \to c \to \sum a$;
			\item $\sigma(x \otimes y) \subseteq \sigma(x).$
		\end{enumerate}
	\end{definition}
	
	We can of course ask the same question as usual: does a universal quasi support datum exist? Indeed it does! 
	
	\begin{theorem}
		Let $\Sp^r(\catK)$ denote the set of \textbf{all} thick right $\otimes$-ideals in $\catK$, and for every $x \in \catK$, let $\sup(x) = \{\catI \in \Sp^r(\catK)\mid x \not\in \calI\}$. We define the topology on $\Sp^r(\catK)$ via setting $\{\sup(x)\}_{x \in \catK}$ as a set of generators for the closed subsets. 
		
		The pair $(\Sp^r(\catK), \sup)$ gives the final quasi support datum on $\catK$ (with supports closed). Explicitly, for any quasi support datum $(Y, \sigma)$, there exists a unique map $f\colon Y \to \Sp^r(\catK)$ given by \[f(y) = \{a \in \catK\mid y \not\in \sigma(a)\}\] such that $\sigma(a) = f\inv(\sup(a)$ for all $a \in \catK$. 
	\end{theorem}
	
	However, we find ourselves in a similar position to that of of Balmer and Ocal in \cite{BO24}, in which they ask the same question for support data of (not-necessarily-tensor) triangulated categories. The universal object $\Sp^r(\catK)$ does indeed classify all thick right $\otimes$-ideals as one may hope, but it does so by itself being the classification! 
	
	\begin{proof}
		First, we must verify that $\Sp^r(\catK)$ is a quasi support datum. Verifications of (a)-(d) are standard. For (e), we have for any $\catI \in \Sp^r(\catK)$ that $x \in \catK$ implies $x \otimes y \in \catK$, hence $\catI \in \sup(x \otimes y) \implies \catI \in \sup(x)$, thus $(\Sp(\catK), \sup)$ is a quasi support datum. Now, given any quasi support datum $(Y, \sigma)$, we must check that $f$ is well-defined, continuous, unique, and it satisfies $\sigma(a) = f\inv(\sup(a)$. This is a familiar drill by now: uniqueness follows by the same argument as in \ref{lem:equalityoffuncs}. To check $f$ is well-defined, it suffices to verify that $f(y)$ is a thick right $\otimes$-ideal. Let $a \in f(y)$, then $y \not\in \sigma(a)$. Condition (e) implies $y \not\in\sigma(a\otimes b)$ for any $b \in \catK$, therefore $a \otimes b \in f(y)$ for any $b \in \catB$, as desired. By definition, we have $f(y) \in \sup(a)$ if and only if $a \not\in f(y)$ if and only if $y \in \sigma(a)$. Therefore, $f\inv(\sup(a)) = \sigma(a)$, as desired, and it follows that $f$ is continuous. 
	\end{proof}
	
	\begin{remark}
		Not to despair, for \cite[Theorem 7.3.1]{NVY22} still provides a way of classifying thick right $\otimes$-ideals (i.e. determining $\Sp(\catK)$) via support methods, although per remark \cite[Remark 7.3.2]{NVY22}, this technique only works when one's quasi-support datum is in fact a noncommutative support datum. Indeed, it seems the lack of well-behavedness of rigidity with one-sided ideals poses a problem for classification. However, the compression of information into the prime right $\otimes$-ideals does not happen for the universal support. Quoting Balmer and Ocal: ``This is not because of whimsical choices in the construction. Being the solution to a universal problem, it is what it is.'' 
		
	\end{remark}

	\section{The complete prime spectra of crossed product categories}\label{sec:crossedprod}
	
	Given a monoidal triangulated category $\catK$ with a group action by a group $G$, one can form the \textit{crossed product category $\catK \rtimes G$}. This construction is an obvious generalization of semidirect products of groups, as well as Benson and Witherspoon's crossed product algebra \cite{BW14} and Plavnik and Witherspoon's modification \cite{PW18}. These crossed product categories are significant for being perhaps the first exhibited examples of monoidal triangulated categories whose (cohomological) support varieties lack the tensor product property. When $\catK$ is a finite tensor category, Bergh, Plavnik and Witherspoon demonstrated in \cite{BPW24} showed that the cohomological support variety of $\catK \rtimes G$ does not have the tensor product property in general. Further, Huang and Vashaw in \cite{HV25} described $\Spc(\catK\rtimes G)$ in terms of $\Spc(\catK)$ when $\catK$ is Noetherian and left or right dualizable, and gave a more general description in terms of $G$-stable ideals. We review this now.
	
	\begin{definition}
		Let $\catK$ be a monoidal triangulated category and $G$ be a group acting on $\catK$, that is, we have a $\otimes$-functor $\text{Cat}(G) \to \Aut_{\otimes, \Delta}(\catK)$ viewing $G$ as a monoidal category with objects the elements of $G$, only identity morphisms, and tensor product encoding group multiplication. For every $x \in \catK$ and $g \in G$, we denote the action of $g$ on $x$ by ${}^gx$. 
		
		The \textit{crossed product category} is the direct sum \[\catK \rtimes G := \bigoplus_{g \in G}\catK\] indexed by elements of $G$ as additive categories. If $x\in \catK$ is an object of $\catK$ then the corresponding object of $\catK \rtimes G$ is denoted $x \boxtimes g$, reflecting the fact that if $\catK$ is $k$-linear, $\catK \rtimes G$ is the Deligne tensor product. We have \[\Hom_{\catK \rtimes G}(x \boxtimes g, y \boxtimes h) := \begin{cases} \Hom_{\catK}(x,y) & \text{if $g = h$} \\ \emptyset & \text{else}\end{cases}\] The monoidal tensor product is defined via \[(x \boxtimes g)\otimes (y \boxtimes h) := (x \otimes {}^gy \boxtimes gh).\] The monoidal unit for $\catK \rtimes G$ is the object $\bbone \boxtimes 1_G$. We view $\catK = \catK \rtimes 1$ as a monoidal triangulated subcategory of $\catK \rtimes G$, which gives a canonical inclusion $\catK \hookrightarrow \catK \rtimes G$. See \cite{HV25} for details. 
	\end{definition}

	\begin{prop}{\cite[Proposition 5.4]{HV25}}
		Let $G$ be a group acting on a monoidal triangulated category $\catK$. Let $\catI$ be a $G$-ideal of $\catK$, i.e. a $G$-invariant thick $\otimes$-ideal and $\catJ$ be a thick $\otimes$-ideal of $\catK\rtimes G$. We have:
		\begin{enumerate}
			\item $\catI \rtimes G$ is a thick $\otimes$-ideal of $\catK \rtimes G$;
			\item $\catJ \cap \catK$ is a $G$-ideal of $\catK$, regarding $\catJ$ as the subcategory $\catJ \boxtimes 1_G\subseteq \catK$;
			\item The maps given by (a) and (b) define a bijection between $G$-ideals of $\catK$ and thick $\otimes$-ideals of $\catK \rtimes G$. Moreover, this bijection sends primes to primes, and restricts to a homeomorphism $\Spc(\catK \rtimes G) \cong G\on{-}\Spc(\catK)$.
		\end{enumerate}
	\end{prop}

	\begin{theorem}\label{thm:crossedcprimes}
		With the setup as above, the homeomorphism $\Spc(\catK \rtimes G) \cong G\on{-}\Spc(\catK)$ restricts to a homeomorphism $\Spc\cp(\catK \rtimes G) \cong G\on{-}\Spc\cp(\catK)$. 
	\end{theorem}
	\begin{proof}
		We must show the maps in \cite[Theorem 5.4]{HV25} preserve complete primality and that the restriction is continuous. First, suppose $\catP \in\Spc\cp(\catK)$ is completely prime, and suppose $(x \boxtimes g)\otimes (y \boxtimes h) \in \catP \rtimes G$. This product is isomorphic to $(x\otimes {}^gy \boxtimes gh)$, and by $G$-invariance, $x \otimes {}^gy \in \catP$. Therefore, either $x$ or ${}^gy$ belongs to $\catP$, hence either $x \boxtimes g$ or ${}^gy \boxtimes h$ belongs to $\catP \rtimes G$. If the former holds, we are done, and if the latter holds, by $G$-invariance we have ${}^gy \boxtimes h \in \catP \rtimes G$ if and only if $y \boxtimes h \in \catP \rtimes G$, so we are done.
		
		Conversely, suppose $\catQ \in \Spc\cp(\catK \rtimes G)$ is completely prime. By \cite[Theorem 5.4]{HV25} we know that $\catQ = \catP \rtimes G$ for some prime $\catP \in \Spc(\catK)$, and we must show $\catP$ is completely prime. Let $x \otimes y \in \catP$, then $(x \boxtimes 1_G)\otimes (y \boxtimes 1_G) = \in \catP \rtimes G$. Therefore either $(x \boxtimes 1_G)$ or $(y \boxtimes 1_G)$ belongs to $\catP \rtimes G$, so either $x$ or $y$ belongs to $\catP$, as desired. Thus the map descends to a bijection between $\Spc\cp(\catK\rtimes G)$ and $G\on{-}\Spc\cp(\catK)$. 
		
		Finally, for continuity, if we denote the bijection $\Phi\colon \Spc\cp(\catK \rtimes G) \to G\on{-}\Spc\cp(\catK),$ with inverse $\Psi$, then it is straightforward to verify that \[\Phi\inv(\supp\cp(x)) = \supp\cp(x\boxtimes g)\] for every $x \in \catK$ and $g \in G$, and similarly, \[\supp\cp(x) = \Psi\inv(\supp\cp(x\boxtimes g)).\] Thus $\Phi$ is a homeomorphism, as desired.
	\end{proof}
	
	We can say some things about one-sided ideals of $\catK \rtimes G$. Here, the story diverges slightly for left and right thick $\otimes$-ideals, and $G$-invariance no longer plays a role. 
	
	\begin{prop}\label{prop:crossedonesideds}
		Let $G$ be a group acting on $\catK$.
		\begin{enumerate}
			\item Let $\catI$ be a (not necessarily $G$-invariant!) thick right $\otimes$-ideal of $\catK$. Then $\catI \rtimes G$ is a thick right $\otimes$-ideal of $\catK \rtimes G$. Conversely, if $\catJ$ is a thick right $\otimes$-ideal of $\catK \rtimes G$, then $\catJ \cap \catK$ is a thick right $\otimes$-ideal of $\catK$. These constructions induce a bijection between the sets of thick right $\otimes$-ideals of $\catK$ and $\catK\rtimes G$.
			\item Let $\catI$ be a thick left $\otimes$ ideal of $\catK$. For each $g \in G$, we define a thick left $\otimes$-ideal $\Delta \catI$ of $\catK$ as follows: the $\catK \boxtimes g$-component of $\catK \rtimes G$ is \[{}^g\catI \boxtimes g \subseteq \catK \boxtimes g.\]
			Conversely, if $\catJ$ is a thick left $\otimes$-ideal of $\catK \rtimes G$, then $\catJ \cap \catK$ is a thick left $\otimes$-ideal of $\catK$. These constructions induce a bijection between the sets of thick left $\otimes$-ideals of $\catK$ and $\catK \rtimes G$.
		\end{enumerate}
	\end{prop}
	\begin{proof}
		Parts (a) and (b) are straightforward verifications, and follow very similarly to those of \cite[Proposition 5.4]{HV25}. The key point is that in (a), $x \otimes g \in \catI \rtimes G$ if and only if $x \otimes gh \in \catI \rtimes G$ for all thick right $\otimes$-ideals $\catI$ of $\catK$, while in (b), $x \otimes g \in \catI' \rtimes G$ if and only if ${}^hg \otimes hg \in \catI' \rtimes G$ for all thick left $\otimes$-ideals $\catI'$ of $\catK$.
	\end{proof}
	
	It is unclear whether a complete analogue of Theorem \ref{thm:crossedcprimes} holds for one sided prime thick $\otimes$-ideals. We now the crossed product construction to exhibit a category with no completely prime two-sided thick $\otimes$-ideals, as well as non-functoriality of both the noncommutative Balmer spectrum and the category of noncommutative support data. 
	
	\begin{theorem}\label{ex:nocps}
		Set $\catK = D_b(k) \oplus D_b(k)$. The crossed product category $\catK \rtimes C_2$, with $C_2$ acting by permuting indices, contains no completely prime ideals. Moreover, the inclusion $F\colon  \catK \hookrightarrow \catK \rtimes C_2$ does not induce a well-defined map $\Spc(F): \Spc(\catK \rtimes C_2) \to \Spc(\catK)$ nor a well-defined functor $\catS(F): \catS(\catK \rtimes C_2) \to \catS(\catK)$. 
	\end{theorem}
	
	\begin{proof}
		We have that $\Spc(\catK)$ is simply two closed points, $\catP_1 = (0, D_b(k))$ and $\catP_2 = (D_b(k), 0)$. \cite[Proposition 7.2]{HV25} implies the $C_2$-prime spectrum of $\Spc(\catK)$ is only one point, the zero ideal of $\catK$. Note the zero ideal is $C_2$-prime but not prime, e.g. $\catP_1 \otimes \catP_2 = 0$. Therefore it is not completely prime either, so the zero ideal of $\catK \rtimes C_2$ is prime but not completely prime, and $\catK \rtimes C_2$ has no complete primes.

		Next, consider the inclusion $F\colon  \catK \hookrightarrow \catK \rtimes C_2$ defined by $x \mapsto x \boxtimes 1_G$. This is a monoidal triangulated functor. We have that $F\inv(\{0\}) = \{0\}$, which is a nonprime thick $\otimes$-ideal, so the map $F\inv$ does not send primes to primes.  
		
		Now, we show $(\Spc(\catK), \supp \circ F)$ is not a support theory for $\catK$. Indeed, let $x_1 = (k[0], 0)$ and $x_2 = (0, k[0])$. We have $x_1 \otimes x_2 = 0$, and since $\catK$ is symmetric, we have $x_1 \otimes \catK \otimes x_2 = \{0\}$. On the other hand, $\supp(x_1 \boxtimes 1_G) = \supp(x_2 \boxtimes 1_G)  = \Spc(\catK \rtimes C_2) = \ast$. Therefore, we have $(\supp\circ F)(x_1) \cap (\supp\circ F)(x_2) = \ast$, but $\bigcup_{x \in \catK}(\supp\circ F)(x_1 \otimes x \otimes x_2) = \emptyset$. Thus $(\Spc, \supp\circ F)$ does not satisfy axiom (e) of the support data axioms, as desired. 
		
	\end{proof}
	
	\begin{remark}\label{ex:nilps}
		The category $\catK \rtimes C_2$ has nilpotent objects - for instance, $x = (k[1], 0) \boxtimes g$ (where $\langle g\rangle = C_2$) satisfies $x \otimes x = 0$. In particular, the zero ideal is prime, therefore is an intersection of non-complete primes, but not radical. 
		
		It is easy to see the above example generalizes to any finite group $G$ and essentially small tensor-triangulated category $\catK$ with $G$ acting on $\catK' := \catK^{\oplus |G|}$ by permuting indices, a categorical ``wreath product'' when $G$ is a symmetric group.
	\end{remark}

	\section{Duality and duo-ness}
	
	For the next two sections, we shift gears towards functoriality and analyze the induced maps on complete prime spectra. Earlier we asked the question of when an induced map is injective or surjective, with injectivity arising from essentially surjective functors. To discuss surjectivity, we restrict our attention to rigid monoidal triangulated categories, which we review now. 
	
	\begin{definition}\label{def:rigid}
		An object $x \in \catK$ is \textit{rigid} if it has left and right \textit{duals}, and $\catK$ is a \textit{rigid} triangulated category if all objects are rigid. Explicitly, there exist exact contravariant functors $(-)^\vee\colon \catK\op\to \catK$ and ${}^\vee(-)\colon  \catK\op \to \catK$, the \textit{left} and \textit{right dual functors} respectively, such that \[(x \otimes - ) \dashv  (x^\vee \otimes -)\] and \[(- \otimes {}^\vee x) \dashv (- \otimes x)\] for all objects $x \in \catK$. These adjunctions come equipped with natural transformations:
		\begin{itemize}
			\item $\on{ev}_x\colon x^\vee \otimes x \to \bbone$ (left evaluation);
			\item $\on{coev}_x\colon \bbone \to x \otimes x^\vee$ (left coevaluation);
			\item $\on{ev}'_x\colon x\otimes {}^\vee x \to \bbone$ (right evaluation);
			\item $\on{coev}_x'\colon \bbone \to {}^\vee x \otimes x$ (right coevaluation).
		\end{itemize}
		which satisfy \[(\id_x \otimes \ev)\circ(\coev_x \otimes \id_x) = \id \text{  and  } (\ev_{x} \otimes \id_{x^\vee})\circ (\id_{x^\vee}\otimes \coev_x)\] for left duals, and \[(\ev_x'\otimes \id_x )\circ (\id_x \otimes \coev_x') = \id \text{  and  } (\id_{{}^\vee x} \otimes \ev_x')\circ (\coev'_x \otimes \id_{{}^\vee x})\] for right duals. Left and right duals are unique up to isomorphism and we have $x \cong ({}^\vee x)^\vee \cong {}^\vee(x^\vee)$. The (co)evaluation identities imply that $x$ is a direct summand of $x \otimes x^\vee \otimes x$ and $ x \otimes {}^\vee x \otimes x$. Via adjunction, one has that $(x \otimes y)^\vee \cong y^\vee \otimes x^\vee$ and ${}^\vee(x \otimes y) \cong {}^\vee y \otimes {}^\vee x$. See \cite[Section 2.10]{EGNO15} for details. 
	\end{definition}
	
	In the presence of rigidity, left and right ideals are dual.
	
	\begin{prop}\label{prop:dualityofideals}
		If $\catK$ is rigid and $\catI$ is a thick right $\otimes$-ideal, then $\catI^\vee$ and ${}^\vee\catI$ are thick left $\otimes$-ideals. 
	\end{prop}
	\begin{proof}
		This easily follows from the identities $(x \otimes y)^\vee \cong y^\vee \otimes x^\vee$ and ${}^\vee(x \otimes y) \cong {}^\vee y \otimes {}^\vee x$. 
	\end{proof}
	
	In the noncommutative setting, rigidity needs to be handled with slightly more care, as left duals and right duals may not coincide. When $x$ is rigid, in the noncommutative setting we no longer have that $x^{\otimes n} = 0$ implies $x = 0$. 
	However, if $\catK$ is rigid, the existence of nilpotent elements in $\catK$ implies that $(\Spc (\catK), \supp)$ no longer satisfies the tensor product property. 
	
	\begin{theorem}{\cite[Theorem 4.2.1]{NVY22b}}
		If $\catK$ is rigid and $(\Spc(\catK), \supp)$ satisfies the tensor product property, then $\catK$ has no nontrivial $\otimes$-nilpotent elements. 
	\end{theorem}
	
	In fact, we require a weaker condition for the theorem to hold: every object of $\catK$ need only have left \textit{or} right duals. The authors in \cite{NVY22b} prove that the tensor product property holds if and only if all thick prime $\otimes$-ideals of $\catK$ are completely prime (\cite[Theorem 3.1.1]{NVY22b}) and show these conditions are satisfied when every thick right $\otimes$-ideal of $\catK$ is, in fact, two-sided.
	
	\begin{theorem}{\cite{NVY22b}}\label{thm:nonilpotents}
		Let $\catK$ be a monoidal triangulated category in which every thick right ideal is two-sided. Then every prime ideal of $\catK$ is completely prime, and as a consequence, $(\Spc(\catK), \supp)$ satisfies the tensor product property. 
	\end{theorem}
	
	We give an easy equivalent formulation of the two-sided condition in the presence of rigidity, also bringing left $\otimes$-ideals into the picture as well. 
	
	\begin{prop}\label{prop:rigidduo}
		Let $\catK$ be rigid. The following are equivalent:
		\begin{enumerate}
			\item Every thick right $\otimes$-ideal of $\catK$ is two-sided;
			\item Every thick right $\otimes$-ideal $\catI$ of $\catK$ satisfies the following property: $x \in \catI$ if and only if $x^\vee \in \catI$;
			\item Every thick right $\otimes$-ideal $\catI$ of $\catK$ satisfies the following property: $x \in \catI$ if and only if ${}^\vee x \in \catI$;
			\item Every thick left $\otimes$-ideal of $\catK$ is two-sided;
			\item Every thick left $\otimes$-ideal $\catI$ of $\catK$ satisfies the following property: $x \in \catI$ if and only if $x^\vee \in \catI$;
			\item Every thick left $\otimes$-ideal $\catI$ of $\catK$ satisfies the following property: $x \in \catI$ if and only if ${}^\vee x \in \catI$. 
		\end{enumerate}
	\end{prop}
	\begin{proof}
		First we show (a) implies (b) and (c). Suppose every thick right $\otimes$-ideal of $\catK$ is two-sided. It suffices to show that if $x \in \catI$, then $x^\vee \in \catI$ and ${}^\vee x \in \catI$, the other implications follow by replacing $x$ with $x^\vee$ or ${}^\vee x$. Since $\catI$ is two-sided, $x^\vee \otimes x \otimes x^\vee \in \catI$. Because $\catI$ is thick and $x^\vee$ is a summand of $x^\vee \otimes x \otimes x^\vee $, $x^\vee \in \catI$. The right dual case follows similarly. (b) implies (c) is obvious. 
		
		Next we show (b) implies (a) suppose $\catI$ is a right $\otimes$-deal that is closed under left duals.  Let $x \in \catI$ and $a \in \catK$. We have $a \otimes x \in \catI$ if and only if $(a \otimes x)^\vee \in \catI$, and $(a\otimes x)^\vee \cong x^\vee \otimes a^\vee \in \catI$, since $x^\vee \in \catI$ by hypothesis. (c) implies (a) follows similarly. The equivalency of (d)-(f) follows from Proposition \ref{prop:dualityofideals}.
	\end{proof}
	
	\begin{definition}
		We say that a monoidal triangulated category $\catK$ satisfying that all thick right (resp. left) $\otimes$-ideals are two sided is \textit{right (resp. left) duo}. If $\catK$ is both right and left duo, we say it is \textit{duo}. In this setting, all prime and completely prime ideals also coincide, allowing more of the ``standard'' primality arguments of tt-geometry. Proposition \ref{prop:rigidduo} asserts that rigid right duo and left duo monoidal triangulated categories in fact coincide. It is an unknown whether there exist monoidal triangulated categories which are right but not left duo, and vice versa. 
	\end{definition}
	
	We record the following obvious yet important fact. 
	
	\begin{prop}
		Suppose $\catK$ is right (resp. left) duo. Then for all $x \in \catK$, $\langle x\rangle_r = \langle x \rangle$ (resp. $\langle x \rangle_l = \langle x\rangle$).
	\end{prop}
	\begin{proof}
		Clearly we have $\langle x\rangle_r \subseteq \langle x\rangle$ since a two-sided ideal is a right ideal. Since $\langle x\rangle_r$ is two-sided, if $\langle x\rangle_r \neq \langle x \rangle$, this would contradict minimality of $\langle x \rangle$ as two-sided $\otimes$-ideal. The left case follows analogously. 
	\end{proof}
	
	We note one general criterion for right duo-ness for well-behaved finite dimensional Hopf algebras. This is used in the computations in \cite[Section 5]{NVY22b}; we record it for completeness. Note that by \cite{NVY24}, one always has a homomorphism into Proj of the categorical center of the cohomology ring, but it is by no means a homeomorphism in general. 
	
	\begin{prop}\label{thm:hopf}
		Let $A$ be a finite-dimensional Hopf algebra over a field $k$ whose Balmer spectrum is homeomorphic to its cohomological support variety, i.e. we have a homeomorphism \[\Spc(\stab(A)) \cong \on{Proj}(\on{H}^\sbull(A, k)).\] If the following conditions hold: 
		\begin{itemize}
			\item The cohomology ring $\on{H}^\sbull(A, k)$ is finitely generated;
			\item For all $M \in \stab(A)$, $\Ext^\sbull(M,M)$ is a finitely generated $\on{H}^\sbull(A, k)$-module;
			\item \cite[Assumption 7.2.1]{NVY22} holds. 
		\end{itemize}
		then $\stab(A)$ is right duo. 
	\end{prop}
	\begin{proof}
		This directly follows from \cite[Theorem 7.4.3]{NVY22} and \cite[Theorem 6.2.1]{NVY22}, since in this case, $(\Spc(\stab(A)), \supp)$ is both a quasi-support datum and support datum which satisfy both the two-sided and right faithfulness and realization properties, as well as \cite[Assumption 7.2.1]{NVY22}. Right duo-ness implies duo-ness since $\stab(A)$ is rigid. 
	\end{proof}
	
	The finite generation conditions are conjectured to hold for all finite tensor categories by Etingof and Ostrik, see \cite{EO04}. The given homeomorphism on Balmer spectra is conjectured to hold for all quasi-triangular Hopf algebras by Nakano, Vashaw, and Yakimov, see \cite{NVY24} (more generally, they predict a modification of this homeomorphism holding for all finite tensor categories). If both conjectures hold, it follows that $\stab(A)$ is duo for any finite dimensional quasi-triangular Hopf algebra $A$. Their conjecture would also provide a positive answer to a question asked by Bergh, Plavnik, and Witherspoon \cite{NVY24} on whether the cohomological support of a braided finite tensor category with finitely generated cohomology has the tensor product property. 
	
	\section{Surjectivity and tensor-nilpotence}\label{sec:surj}
	
	We return to the question of determining when the induced map \[\varphi:=\Spc\cp(F)\colon  \Spc\cp(\catL) \to \Spc\cp(\catK)\] is surjective by adapting the results of \cite{Bal18} to the noncommutative duo setting. 
	
	\begin{definition}\label{def:tensornil}
		A morphism $f\colon x\to y$ is called \textit{$\otimes$-nilpotent} if $f^{\otimes n}\colon x^{\otimes n} \to y^{\otimes n}$ is zero for some $n \geq 1$ (note this can occur even if $\catK$ has no nonzero nilpotent elements!). We say $f\colon x\to y$ is \textit{left (resp. right) $\otimes$-nilpotent on an object $z$} in $\catK$ if there exists $n \geq 1$ such that $f^{\otimes n}\otimes z$ (resp. $z \otimes f^{\otimes n}$) is the zero morphism. If $f$ is both left and right $\otimes$-nilpotent on an object $z$, we simply say $f$ is $\otimes$-nilpotent on $z$. 
		
		In particular, $f$ is \textit{left (resp. right) $\otimes$-nilpotent on its cone} if there exists $n \geq 1$ for which $f^{\otimes n} \otimes \cone(f) = 0$ (resp. $\cone(f) \otimes f^{\otimes n} = 0$).
		
		We also require a weaker notion of $\otimes$-nilpotency on objects for the non-commutative setting. A morphism $f\colon x \to y$ is \textit{weakly $\otimes$-nilpotent on an object $z$} if there exist positive integers $k, i_i, i_2, \dots, i_k \geq 1$ such that $z \otimes f^{i_1} \otimes z \otimes f^{i_2} \otimes \cdots \otimes f^{i_k} \otimes z = 0$. Clearly, left or right $\otimes$-nilpotence on an object implies weak $\otimes$-nilpotence.  
	\end{definition}
	
	\begin{remark}\label{rmk:failureoftriang}
		Let $f\colon x \to y$ be a morphism in $\catK$. Contrary to the classical case, in general it is unclear if the collection of all objects $z \in \catK$ for which $f$ is left (resp. right) $\otimes$-nilpotent on $z$ is a thick $\otimes$-ideal. Though this set is clearly multiplicative on the right (resp. left) and thick, it may not be a triangulated subcategory! The key obstruction is that \cite[Lemma 2.11]{Bal10} shows that this set is a thick $\otimes$-ideal, but uses symmetric structure. However, we do not have any explicit counterexamples of instances when the collection in question is not a triangulated subcategory - finding such a counterexample is a question of interest. 
	\end{remark}
	
	This poses a problem for replicating the results of \cite{Bal18} verbatim, but fear not, the next proposition asserts that such problems do not occur in the situation we are most concerned with: homotopy fibers. 
	
	\begin{prop}\label{prop:tensnilpotence}
		Let $\xi\colon w \to \bbone$ be a morphism in $\catK$ such that $\xi \otimes \cone(\xi) = 0$. Then the cone of $\xi^{\otimes n}$ generates the same thick right $\otimes$-ideal for all $n$:
		\[\langle \cone(\xi) \rangle_r = \{z \in \catK \mid \xi \text{ is left $\otimes$-nilpotent on $z$}\} = \langle \cone(\xi^{\otimes n})\rangle_r.\] In particular, the subcategory $\{z \in \catK \mid \xi \text{ is left $\otimes$-nilpotent on z}\}$ is a thick right $\otimes$-ideal of $\catK$. 
		
		Moreover, if $z \in\catK$ satisfies that $\xi$ is weakly $\otimes$-nilpotent on $z$, then there exists an $n \geq 1$ for which $z^{\otimes n}$ is $\otimes$-nilpotent on $z$. 
	\end{prop}
	\begin{proof}
		The assumption $\xi \otimes \cone(\xi) = 0$ implies that \[\cone(\xi) \in \{z \in \catK \mid \xi \text{ is left $\otimes$-nilpotent on $z$}\}.\] On the other hand, if the morphism $\xi^{\otimes n} \otimes z$ is zero then the exact triangle \[w^{\otimes n} \otimes z \xrightarrow{\xi^{\otimes n} \otimes z = 0} z \to \cone(\xi^{\otimes n})\otimes z \to \Sigma w^{\otimes n} \otimes z\] implies $z$ is a summand of $\cone(\xi^{\otimes n})\otimes z$, thus $z \in \langle \cone(\xi^{\otimes n})\rangle_r.$ Finally, it follows from induction and application of the octahedral axiom that $\cone(f^{\otimes n}) \in \langle \cone(f) \rangle_r$, in particular it fits in an exact triangle \[\cone(f^{\otimes n - 1}) \otimes w \to \cone(f^{\otimes n}) \to \cone(f) \to \Sigma \cone(f^{\otimes n - 1}) \otimes w,\] for all $n \geq 1$, using that $\cone(f) \otimes w \cong \cone(f \otimes w)$.
		Thus $\cone(\xi^{\otimes n}) \in \langle\cone(\xi)\rangle_r$. Thus, we have a chain of inclusions \[\langle\cone(\xi)\rangle_r \subseteq \{z \in \catK \mid \xi \text{ is left $\otimes$-nilpotent on $z$}\} \subseteq \bigcup_{n \geq 1} \langle \cone(\xi^{\otimes n})\rangle_r \subseteq \langle\cone(\xi)\rangle_r.\] This proves the first claim - in particular, the middle set is indeed a thick triangulated right $\otimes$-ideal. 
		
		For the final statement, if $\catK$ is duo, because $\langle \cone(\xi)\rangle_l$ and $\langle \cone(\xi)\rangle_r$ coincide, left and right $\otimes$-nilpotence coincide for $\xi$, and $\langle \cone(\xi) \rangle = \{z \in \catK \mid \xi \text{ is $\otimes$-nilpotent on $z$}\}$. Now, we have a similar exact triangle, where we've multiplied by $\xi$ on the left: \[w\otimes z \otimes w^{i_1} \otimes \cdots \otimes w^{i_k}\otimes z \xrightarrow{\xi \otimes z \otimes \xi^{i_1} \otimes \cdots \otimes \xi^{i_k} \otimes z =0} z^{\otimes (k+1)} \to \cone(\xi \otimes z \otimes \xi^{i_1} \otimes \cdots \otimes \xi^{i_k} \otimes z) \to \dots, \] therefore $z^{\otimes k+1}$ is a direct summand of $\cone(\xi \otimes z \otimes \xi^{i_1} \otimes \cdots \otimes \xi^{i_k} \otimes z)$ and by a similar argument as before, belongs to $\langle\cone(\xi)\rangle_r$. Multiplying by $\xi$ on the right shows $z^{\otimes k+1}$ belongs to $\langle\cone(\xi)\rangle_l$.
	\end{proof}
	
	As a result, we obtain a key corollary allowing us to translate questions about right $\otimes$-ideals into questions about left $\otimes$-nilpotence in the presence of rigidity via our well-behaved homotopy fibers. This is the only lemma explicitly requiring duo-ness; whether the condition may be removed would certainly be interesting to know. 
	
	\begin{corollary}\label{cor:tensidealsnilpotents}
		Let $x \in \catK$ be a rigid object in a (not necessarily rigid) monoidal triangulated category $\catK$. Choose $\xi_x$ a homotopy fiber of the coevaluation morphism $\coev_x$, i.e. choose an exact triangle \[w_x \xrightarrow{\xi_x} \bbone \xrightarrow{\coev_x} x \otimes x^\vee \to \Sigma w_x.\] Then the thick right $\otimes$-ideal $\langle x \rangle_r$ generated by $x$ is exactly the subcategory on which $\xi_x$ is $\otimes$-nilpotent on the left, that is, \[\langle x \rangle_r = \{z \in \catK \mid \xi_x^{\otimes n} \otimes z = 0 \text{ for some } n \geq 1\}.\] Moreover, for every $n \geq 1$, the morphism $\xi_x^{\otimes n}$ is left $\otimes$-nilpotent on its cone. 
	\end{corollary}
	\begin{proof}
		Consider the exact triangle obtained by tensoring the above with $x$ on the right: \[ w_x \otimes x\xrightarrow{\xi_x \otimes x} x \xrightarrow{\coev_x \otimes x} x \otimes x^\vee \otimes x \to \Sigma w_x \otimes x.\] By properties of left duals, $\coev_x \otimes x$ is a (split) monomorphism, which forces $\xi_x \otimes x = 0$. Therefore, $\xi_x \otimes \cone(\xi_x) \cong  \xi_x \otimes x \otimes x^\vee = 0$. Therefore, we can apply Proposition \ref{prop:tensnilpotence} to $\xi = \xi_x$. This gives us that $\langle \cone(\xi_x)\rangle_r = \{z \in \catK \mid \xi_x \text{ is left $\otimes$-nilpotent on $z$}\}.$ However, $\langle\cone(\xi_x)\rangle_r = \langle x \otimes x^\vee\rangle_r$ and since $x$ is a direct summand of $x \otimes x^\vee \otimes x$, it easily follows that  $\langle x \rangle_r = \langle x \otimes x^\vee \rangle_r$. Thus $\langle x\rangle_r = \{z \in \catK \mid \xi_x \text{ is left $\otimes$-nilpotent on $z$}\}$ as desired. The final statement also follows directly from Proposition \ref{prop:tensnilpotence}. 
	\end{proof}

	\begin{remark}\label{rmk:nilpotentflipsides}
		One may interchange ``left'' and ``right'' in the statement of Corollary \ref{cor:tensidealsnilpotents} and its proof. Therefore, for every $n \geq 1$, the morphism $\xi_x^{\otimes n}$ is $\otimes$-nilpotent on its cone. Moreover, by the final statement of Proposition \ref{prop:tensnilpotence}, if $\xi_x$ is weakly $\otimes$-nilpotent on $z$, there exists an $n \geq 1$ for which $z^{\otimes n}$ belongs to $\langle x \rangle_r$ and $\langle x \rangle_l$. This observation will be critical in the sequel! 
		
		It follows that if $\catK$ is duo, then for all $x \in \catK$, $z \in \catK$ is left $\otimes$-nilpotent on $\xi_x$ if and only if $z$ is right $\otimes$-nilpotent on $\xi_x$. 
	\end{remark}
	
	We require one final lemma, a surjectivity argument relying on the previous corollary, Proposition \ref{prop:primelift}, and duo-ness. 
	
	\begin{lemma}\label{lem:technical}
		Suppose $\catK$ is rigid, $\catL$ is duo, and let $F\colon \catK \to \catL$ be a monoidal triangulated  functor. Choose for every $x \in \catK$ an exact triangle as in Corollary \ref{cor:tensidealsnilpotents}, i.e. a homotopy fiber of the coevaluation morphism $\coev_x$. Let $\catP \in \Spc\cp(\catK)$ be a complete prime. Suppose that $\catP$ satisfies the following technical condition: \[\text{For all $x \in \catP$, all $s \in \catK \setminus \catP$, and all $n \geq 1$, we have $F(\xi_x^{\otimes n}\otimes s )\neq 0$ }.\] Then $\catP$ belongs to the image of $\varphi\colon \Spc\cp(\catL) = \Spc(\catL) \to \Spc\cp(\catK)$. 
	\end{lemma}
	\begin{proof}
		Let $\catI \subseteq \catL$ be the right $\otimes$-ideal generated by $F(\catP)$. We claim that $\catI$ is equivalently \[\catI' := \{y \in \catL \mid \exists \,x \in \catP \text{ such that } y \in \langle F(x)\rangle_r\}.\] Indeed, it is clear that $F(\catP) \subseteq \catI' \subseteq \catI$ straight from the definitions. It suffices to show that $\catI'$ is a thick right $\otimes$-ideal. It is clearly thick and a right $\otimes$-ideal. For closure under cones, suppose $y_1 \to y_2 \to y_3 \to \Sigma y_1$ is exact in $\catL$ and $y_i \in \langle F(x_i)\rangle_r$ for $x_i \in \catP$ and $i \in \{1,2\}$. Then \[y_3 \in \langle y_1, y_2\rangle_r \subset \langle F(x_1), F(x_2)\rangle_r = \langle F(x_1 \oplus x_2)\rangle_r,\] and since $x_1 \oplus x_2 \in \catP$, $y_3 \in \catI'$. We conclude \[\catI = \langle F(\catP)\rangle_r =\{y \in \catL \mid \exists \,x \in \catP \text{ such that } y \in \langle F(x)\rangle_r\}.\]
		
		Now, for every $x \in \catK$, $F$ sends the exact triangle of Corollary \ref{cor:tensidealsnilpotents} to an exact triangle in $\catL$: \[F(w_x) \xrightarrow{F(\xi_x)} \bbone \xrightarrow{\eta_{F(x)}} F(x) \otimes F(x)^\vee \to \Sigma F(w_x).\] Note that $F$ preserves duals, as it is a $\otimes$-functor. Corollary \ref{cor:tensidealsnilpotents} implies \[ \langle F(x)\rangle_r = \{y \in \catL \mid F(\xi_x)^{\otimes n} \otimes y = 0 \text{ for some } n \geq 1\}.\] Combining this with the description of $\catI = \langle F(\catP)\rangle_r$, we have \[\langle F(\catP)\rangle_r = \{y\in \catL \mid F(\xi_x)^{\otimes n} \otimes y = 0 \text{ for some } n \geq 1 \text{ and some } x \in \catP\}.\] Recall that we have assumed $\catL$ is duo, hence $\langle F(\catP)\rangle_r = \langle F(\catP)\rangle$. 
		We set $S := \catK \setminus \catP$; it follows that if $s \in S$, then $F(s)$ cannot belong to $\catI = \langle F(\catP)\rangle$. Indeed, if $F(s) \in \langle F(\catP)\rangle$, then by the above, there exists an $x \in \catP$ and $n \geq 1$ such that \[0 = F(\xi_x)^{\otimes n} \otimes F(s) \cong F(\xi_x^{\otimes n} \otimes s).\] But this directly contradicts our technical assumption on $\catP$. 
		
		Therefore, we have shown that the two-sided $\otimes$-multiplicative class $F(S) = F(\catK \setminus \catP)$ does not meet the $\otimes$-ideal $\catI = \langle F(P)\rangle$ in $\catL$. By the prime lifting trick Proposition \ref{prop:primelift}, there exists a prime $\catQ$ satisfying that $\catI \subseteq \catQ$ and $F(S) \cap \catQ = \emptyset$, and $\catQ$ is completely prime since duo-ness implies all primes are completely prime. These relations mean respectively $\catP \subseteq F\inv(\catQ)$ and $F\inv(\catQ) \subseteq \catP$, hence $\catP = F\inv(\catQ) = \varphi(\catQ)$, as desired.
	\end{proof}
	
	The following theorem is a weaker form of \cite[Theorem 1.3]{Bal18}. Note that $\otimes$-nilpotence of morphisms is detected for instance by conservative functors on compactly generated categories restricted to their compact parts, see \cite[Proposition 7.1]{BG23} for instance. 
	
	\begin{theorem}\label{thm:thm2}
		Assume $\catK$ is rigid and $\catL$ is duo. Suppose that the exact $\otimes$-functor $F\colon \catK \to \catL$ detects $\otimes$-nilpotents of morphisms, i.e. every $f\colon x\to y$ in $\catK$ such that $F(f) = 0$ satisfies $f^{\otimes n} = 0$ for some $n \geq 1$. Then the induced map $\varphi\colon \Spc\cp(\catL)\to \Spc\cp(\catK)$ is surjective.
	\end{theorem}
	\begin{proof}
		
		Suppose that $f\colon \catK \to \catL$ detects $\otimes$-nilpotence of morphisms. Let $\catP \in \Spc\cp(\catK)$ be a complete prime. We show that the technical condition of Lemma \ref{lem:technical} is satisfied, i.e. \[\text{For all $x \in \catP$, all $s \in \catK \setminus \catP$, and all $n \geq 1$, we have $F(\xi_x^{\otimes n}\otimes s )\neq 0$}.\] Let $g = \xi_x^{\otimes n} \otimes s$ be the morphism in the condition for some objects $x \in \catP$ and $s \in \catK \setminus\catP$ and some $n \geq 1$, and suppose for contradiction that $F(g) = 0$. Then by assumption, $g = \xi_x^{\otimes n} \otimes s$ is $\otimes$-nilpotent. In other words, $\xi_x$ is weakly $\otimes$-nilpotent on $s$, therefore by Proposition \ref{prop:tensnilpotence}, $\xi_x^{\otimes k}$ is $\otimes$-nilpotent on $s$ for some $k \geq 1.$ This implies that $s^{\otimes k}$ belongs to $\langle x \rangle$ by Corollary \ref{cor:tensidealsnilpotents}, but $\langle x \rangle \subseteq \catP$ by assumption, and since $\catP$ is completely prime, $s \in \catP$, which directly contradicts the choice of $s \in \catK \setminus \catP$. Thus, we have verified the property of Lemma \ref{lem:technical} holds for the prime $\catP$, therefore $\catP$ belongs to the image of $\varphi$, as desired.

	\end{proof}

	\begin{remark}\label{rmk:cantget14}
		There are two key points which prevent us from fully extending \cite[Theorem 1.4]{Bal18}, a stronger statement that characterizes exactly when $\varphi:\Spc(\catL)\to\Spc(\catK)$ is surjective in generality. In one direction (the direction (a) $\implies$ (b) in the proof), it requires that the collection $\{z \in \catK \mid f \text{ is $\otimes$-nilpotent on }z\}$ be a thick triangulated $\otimes$-ideal, which as noted in \ref{rmk:failureoftriang}, may not be a triangulated subcategory. Additionally, in the (b) $\implies$ (a) direction, it is deduced from the fact that $\xi_x^{\otimes n}$ is $\otimes$-nilpotent on its cone that $g = \xi_x^{\otimes n} \otimes s$ is $\otimes$-nilpotent on its cone as well; this does not hold in our case, even if we replace ``$\otimes$-nilpotent'' with ``weakly $\otimes$-nilpotent.'' 
	\end{remark}

	Now we enter the hall of the right adjoint. We assume $\catK$ is rigid and $f\colon \catK \to \catL$ has a right adjoint $U\colon \catL \to \catK$ (which is not necessarily a $\otimes$-functor). It follows by abstract nonsense that $U$ satisfies the projection formula \[\pi\colon x \otimes U(y) \cong U(F(x)\otimes y)\] for all $x \in \catK$ and $y \in \catL$. In fact, for the projection formula to be an isomorphism, we only require that every object of $\catK$ has a left or right dual. 
	
	The next two proofs follow nearly verbatim from the proofs of \cite[Theorems 1.7 and 1.6]{Bal18} respectively.
	\begin{theorem}\label{thm:thm3}
		Suppose $\catK$ is rigid, $\catL$ is duo, and $F\colon \catK \to \catL$ admits a right adjoint $U\colon \catL \to \catK$. Then the image of the map $\varphi\colon \Spc(\catL) \to \Spc\cp(\catK)$ is exactly the CP-support of the image of the $\otimes$-unit $U(\bbone)\colon $ \[\im(\varphi) = \supp_{\catK}\cp(U(\bbone)).\]
	\end{theorem}
	\begin{proof}
		Let $\catP \in \Spc\cp(\catK)$ be a complete prime. We need to show that $\catP \in \im(\varphi)$ if and only if $\catP \in \supp_\catK\cp(U(\bbone_\catL))$, or equivalently $U(\bbone_\catL) \not\in \catP$. 
		
		First, suppose $\catP = \varphi(\catQ)$ for some $\catQ \in \Spc(\catL)$. Then $\catP = F\inv(\catQ)$ by definition. To show $U(\bbone_\catL) \not\in\catP$, it therefore suffices to show $FU(\bbone_\catL)\not\in\catQ$. But this follows from the unit-counit relation for $F$ and $U$: the object $FU(\bbone_\catL) \cong FUF(\bbone_\catK)$ admits $F(\bbone_\catK) \cong \bbone_\catL$ as a direct summand, and $\bbone_\catL$ belongs to no prime in $\Spc(\catL).$
		
		Now, let $\catP \in \supp_\catK\cp(U(\bbone_\catL))$, or equivalently, $U(\bbone_\catL) \not\in \catP$. We show that $\catP$ satisfies the technical condition of Lemma \ref{lem:technical}, \[\text{for all $x \in \catP$, all $s \in \catK \setminus \catP$, and all $n \geq 1$, we have $F(\xi_x^{\otimes n}\otimes s )\neq 0$}.\] Let $x \in \calP$ and $s \in \catK \setminus \catP$, set $g = \xi_x^{\otimes n} \otimes s$ for some $n \geq 1$, and suppose for sake of contradiction that $F(g) = 0$. By the projection formula with $y = \bbone_\catL$, $UF(g)= U(0) = 0$ implies that $g \otimes U(\bbone_\catL) = 0$. Therefore, we have an exact triangle \[w_x^{\otimes n} \otimes s \otimes U(\bbone_\catL) \xrightarrow{g \otimes U(\bbone_\catL) = 0} s \otimes U(\bbone_\catL) \to \cone(g) \otimes U(\bbone_\catL) \to \Sigma w_x^{\otimes n} \otimes s \otimes U(\bbone_\catL)\] in $\catK$. Since $g\otimes U(\bbone_\catL) = 0$, $s \otimes U(\bbone_\catL)$ is a direct summand of $\cone(g) \otimes U(\bbone_\catL) \in \langle \cone(g) \rangle$. But $\cone(g) = \cone(\xi_{x}^{\otimes n})\otimes s \in \langle\cone( \xi_x^{\otimes n})\rangle$, and by Corollary \ref{cor:tensidealsnilpotents}, $\cone(\xi_x^{\otimes n}) \in \langle x \rangle \subseteq \catP$. Thus, $s \otimes U(\bbone_\catL) \in \catP$. By complete primality, this forces either $s \in \catP$ or $U(\bbone_\catL) \in \catP$, which are both absurd. Thus $F(g) \neq 0$, and we conclude by Lemma \ref{lem:technical} that $\catP \in \im(\varphi)$ as desired.
	\end{proof}
	
	\begin{theorem}\label{thm:thm4}
		Suppose $\catK$ is rigid, $\catL$ is duo, and the exact $\otimes$-functor $F\colon \catK \to \catL$ admits a right adjoint $U\colon \catL \to \catK$. Then $\varphi\colon \Spc(\catL) \to \Spc\cp(\catK)$ is surjective if and only if the monoidal $\otimes$-functor $F\colon \catK \to \catL$ detects $\otimes$-nilpotence of morphisms. 
	\end{theorem}
	\begin{proof}
		In light of Theorem \ref{thm:thm3}, it suffices to show that $f\colon \catK\to\catL$ detects $\otimes$-nilpotence if and only if $\supp_\catK\cp(U(\bbone_\catL)) = \Spc(\catK)$, or equivalently $\langle U(\bbone_\catL)\rangle = \catK$. Since $U$ is lax monoidal, $A := U(\bbone_\catL)$ is a ring object. Let \[J \xrightarrow{\xi} \bbone_\catK \xrightarrow{u} A \to \Sigma J\] be an exact triangle over the unit $u\colon \bbone_{\catK} \to A$, i.e. the unit of the $F \dashv U$ adjunction evaluated at $\bbone_{\catK}$. We have $A \otimes \xi = 0$ since $A\otimes u$ is a split monomorphism retracted by multiplication $A \otimes A \to A$. 
		
		Now a morphism $f\colon x\to y$ satisfies $F(f) = 0$ if and only the composition $x\xrightarrow{f} y \xrightarrow{u \otimes y} A \otimes y$ is zero. Indeed, this follows from the projection formula $A \otimes - \cong UF(-)$ with $y = \bbone_\catL$. Now, this is equivalent to the morphism $f\colon x\to y$ factoring via $\xi \otimes y\colon J \otimes y \to y$, via the exact triangle \[J \otimes y \xrightarrow{\xi \otimes y} y \xrightarrow{u \otimes y} A \otimes y \to \Sigma J \otimes y.\] Since $f\colon x\to y$ factors through $\xi \otimes y$ if and only if $F(f) = 0$, it suffices to show $\xi\colon J \to \bbone_\catK$ is $\otimes$-nilpotent if and only if $\langle A\rangle = \catK$. However, this follows immediately from Proposition \ref{prop:tensnilpotence}, which says that $\langle A \rangle = \{z \in \catK \mid \xi \text{ is $\otimes$-nilpotent on $z$.}\}$, and of course, $\bbone_\catK \in \langle A\rangle$ if and only if $\xi$ is $\otimes$-nilpotent on $\bbone_\catK$, i.e. $\otimes$-nilpotent. 
	\end{proof}
	
	\section{Finale: the universal functorial support theory}
	
	Throughout, we have concerned ourselves with the question of functorial support theories for monoidal triangulated categories. It is surely worth asking then if a universal support theory extending the Balmer spectrum for braided monoidal triangulated categories exists. In the final section, we show the answer is yes, following the blueprint laid out by Reyes \cite{Rey12}. For brevity, we say two objects $x, y\in \catK$ \textit{commute} if there exists an isomorphism $x \otimes y \cong y\otimes x$. 
	
	\begin{definition}
		We say an exact triangle  $x \xrightarrow{f} y \xrightarrow{g} z \xrightarrow{h} \Sigma x$ in $\catK$ \textit{commutes} if it is contained in a monoidal triangulated subcategory of $\catK$ (not necessarily full or thick) with the same tensor unit and a braided structure. If this holds, one has an isomorphism $\gamma: x\otimes y \cong y \otimes x$ satisfying $(\id \otimes f)\circ \gamma = f\otimes \id$ and $f\otimes \id = \gamma \circ (\id \otimes f)$. 
	\end{definition}
	
	\begin{definition}\label{def:partialprime}
		A full subcategory $\catI$ of $\catK$ is a \textit{partial $\otimes$-ideal} if the following hold.
		\begin{enumerate}
			\item $\catI$ is closed under shifts. 
			\item If $x, y\in \catK$ commute, then $x, y\in \catI$ if and only if $x \oplus y \in \catI$.
			\item If $x, y \in\catK$ commute and $x \in \catI$, then $x\otimes y \cong y \otimes x \in\catI$. 
			\item Given a commuting exact triangle $x \xrightarrow{f} y \xrightarrow{g} z \xrightarrow{h} \Sigma x$ in $\catK$, if $x, y\in \catI$, then $z \in\catI$. 
		\end{enumerate}
		
		A proper partial $\otimes$-ideal $\catI$ is \textit{prime} if furthermore, the following holds for all commuting $x, y \in\catK$: if $x\otimes y \in \catI$, then $x \in \catI$ or $y \in \catI$. In particular, one has that for any pair of commuting $x,y \in \catK$, $x \otimes y \in \catI$ if and only if $x \in \catI$ or $y \in \catI$.  
	\end{definition}
	
	We remark that $\catI$ may not be either a triangulated subcategory nor a thick $\otimes$-ideal, and possibly not even additive. For instance, the collection of all nilpotent objects in $\catK$ forms a partial $\otimes$-ideal, but may not be a thick $\otimes$-ideal in noncommutative settings. 
	
	\begin{definition}
		Let $p\Spc(\catK)$ denote the set of all prime partial $\otimes$-ideal of $\catK$. Given any $x \in \catK$, let $p\supp(x) = \{\catP \in p\Spc(\catK) \mid x \not\in \catP\}$. As usual, the compliments of all $p\supp(x)$ generate the Zariski topology of $p\Spc(\catK)$. 
	\end{definition}
	
	\begin{prop}\label{prop:partialidealfunctorial}
		Let $F\colon  \catK \to \catL$ be an exact $\otimes$-functor. If $\catI$ is a partial $\otimes$-ideal of $\catL$, then $F\inv(\catI)$ is a partial $\otimes$-ideal of $\catK$. Moreover, if $\catI$ is prime, then so is $F\inv(\catI)$. 
	\end{prop}
	
	\begin{proof}
		First, note $F\inv(\catI)$ is again full. We evaluate conditions (a)-(d). Condition (a) is straightforward. For condition (b), let $x, y \in \catK$ commute. Then $F(x), F(y) \in \catI$ commute, therefore we have $F(x), F(y) \in \catI$ if and only if $F(x) \oplus F(y) \cong F(x \oplus y) \in \catI$, implying (b). For (c), let $x, y\in \catK$ commute and $x \in F\inv(\catI)$. Then $F(x), F(y) \in \catI$ and $F(x) \in \catI$, hence $F(y \otimes x) \in \catI$, whence (c). Finally for (d), suppose $x \xrightarrow{f} y \xrightarrow{g} z \xrightarrow{h} \Sigma x$ is commuting in $\catK$, with it contained in a braided monoidal subcategory $\catK' \subseteq \catK$. Then its image under $F$ is commuting in $\catL$, since it is contained in $F(\catK')$. Therefore, if $F(x), F(y) \in \catI$, then $F(z) \in \catI$, whence (d).
		
		Finally, we evaluate primality. If $x, y\in \catK$ are commuting with $x \otimes y \in F\inv(\catI)$, then $F(x) , F(y)$ are commuting with $F(x) \otimes F(y) \in \catI$. Therefore, $F(x)$ or $F(y) \in \catI$, so either $x$ or $y \in \catI$, as desired.  
	\end{proof}
	
	\begin{corollary}
		Given a monoidal $\otimes$-functor $f\colon \catK \to \catL$, the map \[\varphi := p\Spc(F)\colon  \Spc^{cp}(\catL) \to p\Spc(\catK),\quad \catQ \mapsto F\inv(\catQ)\] is well-defined, continuous, and for all objects $x \in \catK$, we have \[\varphi\inv\big(p\supp_{\catK}(x)\big) = p\supp_\catL(F(x)).\] This defines a contravariant functor $p\Spc$ from the category of essentially small monoidal triangulated categories to the category of topological spaces. 
	\end{corollary}
	\begin{proof}
		This verification follows analogously to Theorem \ref{thm:functorial} using Proposition \ref{prop:partialidealfunctorial}.
	\end{proof}
	
	The following observation is critical. 
	
	\begin{prop}\label{prop:commsubcats}
		Let $x, y\in\catK$ be (possibly the same) commuting elements of $\catK$. There exists a braided monoidal triangulated subcategory of $\catK$ containing $x$ and $y$. 
	\end{prop}
	\begin{proof}
		Take the subcategory of $\catK$ given by finite coproducts and tensor products of $x$ and $y$ containing only isomorphisms and zero maps for morphisms. This category is trivially triangulated and monoidal, and a braiding can be endowed by choosing an isomorphism $x \otimes y \cong y \otimes x$ and extending it to all products of $x$ and $y$ - this satisfies the hexagon identities by construction. 
	\end{proof}
	
	These subcategories are rather ad-hoc - the point is only that they exist and that their inclusion is an exact $\otimes$-functor, since we require that the subcategory shares the same tensor unit. We now give an important reformulation of partial primeness. 
	
	\begin{prop}
		Let $\catP$ be a full subcategory of $\catK$. The following are equivalent: $\catP$ is a partial prime $\otimes$-ideal of $\catK$ if and only if for every (not necessarily full or thick) braided monoidal triangulated subcategory $\catC$ of $\catK$ with the same tensor unit, $\catP \cap \catC$ is a prime thick $\otimes$-ideal of $\catC$. 
	\end{prop}
	\begin{proof}
		First, let us suppose $\catP$ is a partial prime $\otimes$-ideal of $\catK$ and fix $\catC$ a braided monoidal triangulated subcategory $\catC\subseteq \catK$. From the axioms of partial $\otimes$-ideals, $\catP\cap \catC$ is necessarily triangulated and a $\otimes$-ideal, since the axioms hold for all elements of $\catC$. Note $\catC \cap \catP$  is a thick, full subcategory by assumption. Similarly, it follows by definition that $\catP\cap \catC$ is a prime of $\catC$.
		
		Now let us assume the converse: for every braided monoidal triangulated subcategory $\catC$, that $\catP\cap \catC$ is a prime thick $\otimes$-ideal of $\catC$. We must verify all the conditions of a thick prime $\otimes$-ideal hold. They key point here is Proposition \ref{prop:commsubcats}; existence of triangulated subcategories of $\catK$ containing commuting elements or a commuting exact triangle verifies conditions (a), (b), (c), condition (d) follows in a similar fashion, as does primality. 
	\end{proof}
	
	For instance, the partial $\otimes$-ideal consisting of all nilpotent elements of $\catK$ is contained in all prime partial $\otimes$-ideals. 
	
	\begin{definition}
		Let $\Br(\catK)$ denote the set of all braided monoidal triangulated subcategories of $\catK$ with the same tensor unit. In fact, $\Br(\catK)$ is a poset under inclusion, and because inclusion of monoidal triangulated subcategories is an exact $\otimes$-functor, $\Br(\catK)$ is a sub-poset of $\mathbf{Bmon}_\Delta$, the category of essentially small braided monoidal triangulated categories. 
	\end{definition}
	
	\begin{prop}\label{prop:specifyingpartialprime}
		The following data determines a unique partial $\otimes$-ideal of $\catK$: for every $\catC \in \Br(\catK)$, a prime thick $\otimes$-ideal $\catP_\catC$ of $\catC$, such that for all pairs $\catC_1, \catC_2 \in \Br(\catK)$, $\catP_{\catC_1} \cap \catC_2 = \catC_1 \cap \catP_{\catC_2}$. 
	\end{prop}
	\begin{proof}
		We construct a bijection between the set of coherent assignments $\{\catP_\catC\in \Spc(\catC)\}_{\catC \in \Br(\catK)}$ and the set of prime partial $\otimes$-ideals of $\catK$. Given a prime partial $\otimes$-ideal $\catP$, $\{\catP \cap \catC\}_{\catC \in \Br(\catK)}$ forms a coherent assignment of prime thick $\otimes$-ideals. Conversely, given a set of coherent assignments $\{\catP_\catC\in \Spc(\catC)\}_{\catC \in \Br(\catK)}$, the full subcategory with objects given by $\bigcup_{\catC \in \Br(\catK)} \catC_\catP$ is a full prime partial $\otimes$-ideal. It is straightforward to verify these constructions are inverse to each other, with a key point being that every object of $\catK$ belongs to a subcategory $\catC \in \Br(\catK)$. 
	\end{proof}
	
	The next theorem follows essentially immediately from Proposition \ref{prop:specifyingpartialprime}.
	\begin{theorem}
		We have an functorial isomorphism of topological spaces \[p\Spc(\catK) \cong \varprojlim_{\catC \in \Br(\catK)\op} \Spc(\catC).\] This isomorphism preserves the isomorphism of functors $p\Spc \cong \Spc$ when restricted to the full subcategory $\mathbf{Bmon}_\Delta \subset \mathbf{mon}_\Delta$ of essentially small braided triangulated categories. 
	\end{theorem}
	
	\begin{theorem}
		The functor $p\Spc: \mathbf{mon}_\Delta\op \to \mathbf{Top}$, along with the identity natural transformation $p\Spc \cong \Spc$ upon restriction to $\mathbf{Bmon}_\Delta$ is the right Kan extension of the functor $\Spc: \mathbf{Bmon}_\Delta\op \to \mathbf{Top}$ along the embedding $\mathbf{Bmon}_\Delta \hookrightarrow \mathbf{mon}_\Delta$. 
	\end{theorem}
	\begin{proof}
		The setup is as follows:
		\begin{figure}[H]
			\centering
			\begin{tikzcd}
				\mathbf{Bmon}_\Delta\op \ar[d, hookrightarrow] \ar[r, "\Spc"] & \mathbf{Top}\\
				\mathbf{mon}_\Delta \ar[ur, dotted, "p\Spc"']
			\end{tikzcd}
		\end{figure}
		Let $F\colon  \mathbf{mon}_\Delta\op \to \mathbf{Top}$ be a functor for which there exists a natural transformation $\eta: F \to \Spc$ upon restriction to $\mathbf{Bmon}_\Delta$. We must show there is a unique natural transformation $\delta: F \to p\Spc$ that when restricted to $\mathbf{Bmon}_\Delta$ induces $\eta$. For every essentially small monoidal triangulated category $\catK$ and every braided monoidal triangulated subcategory $\catC \subseteq \catK$, the inclusion induces a continuous map $F(\catK) \to F(\catC)$, and we have a morphism $\eta_\catC: F(\catC) \to \Spc(\catC)$, hence a morphism $F(\catK) \to \Spc(\catC)$. Therefore, we have a cone over the diagram $\Br(\catK)\op$, so by universality of the inverse limit, we have a unique morphism \[\delta_\catK: F(\catK) \to \varprojlim_{\catC \in \Br(\catK)\op} \Spc(\catC) \cong p\Spc(\catK).\] This forms the data of our desired unique natural transformation $\delta: F \to p\Spc$, and by construction, $\delta = \eta$ when restricted to $\mathbf{Bmon}_\Delta$. 
	\end{proof}
	
	In particular, $p\Spc$ is not only universal among functors $\mathbf{mon}_\Delta\op \to \mathbf{Top}$ isomorphic to $\Spc$ upon restriction to braided monoidal triangulated categories, it is universal among functors for which there exists only a natural transformation to $\Spc$. 
	
	\begin{corollary}\label{cor:universalfunctorialfunctor}
		The functor $p\Spc: \mathbf{mon}_\Delta\op \to \mathbf{Top}$ is the final object in the full subcategory of $\on{Fun}(\mathbf{mon}_\Delta\op, \mathbf{Top})$ consisting of functors $F\colon  \mathbf{mon}_\Delta\op \to\mathbf{Top}$ which have a natural transformation to $\Spc$ upon restriction to $\mathbf{Bmon}_\Delta\op$. In particular, it is the final object in the full subcategory consisting of functors whose natural transformations are natural isomorphisms as well, i.e. the fiber category $\iota\inv(\Spc)$, where $\iota$ denotes the inclusion $\mathbf{Bmon}_\Delta\op \hookrightarrow \mathbf{mon}_\Delta\op$.
	\end{corollary}
	
	\begin{remark}
		One may carry out the same process of creating the partial prime spectrum with the phrase ``symmetric'' replacing ``braided,'' but with the additional imposed condition that for $x,y$ to ``commute,'' there must exist an isomorphism $\tau: x \otimes y \cong y \otimes x$ such that $\tau^2 = \id$. 
	\end{remark}

	Now that such a universal spectral theory for essentially small monoidal triangulated categories exists, one may ask if it is universal in the other sense, in that it is universal among support theories satisfying a certain set of support axioms. If it is nonempty, it is, in a rather contrived manner. 
	
	\begin{definition}\label{def:commsupp}
		A \textit{commutative support datum} for a monoidal triangulated category $\catK$ is a map \[\sigma: \catK \to \calX_{cl}(Y)\] for a topological space $Y$ such that the following hold:
		\begin{enumerate}
			\item $\sigma(0) = \emptyset$ and $\sigma(\bbone) = Y$;
			\item $\sigma(x \oplus y) = \sigma(x) \cup \sigma(y)$ for all commuting $x,y\in\catK$;
			\item $\sigma(\Sigma x) = \sigma(x)$;
			\item $\sigma(x) \subseteq \sigma(y) \cup \sigma(z)$ for all commuting triangles $x \to y\to z\to \Sigma x$ in $\catK$;
			\item $\sigma(x \otimes y) = \sigma(x) \cap \sigma(y)$ for all commuting $x,y \in \catK$
		\end{enumerate}
		A \textit{morphism of commutative support data} is a continuous map $f:(Y,\sigma) \to (Z,\tau)$ such that $f\inv(\tau(x)) = \sigma(x)$ for all $x \in \catK$. Let $\catS^{cm}(\catK)$ denote the category of commutative support data on $\catK$. 
	\end{definition}
	
	\begin{theorem}\label{thm:finalcommsuppdata}
		If $p\Spc(\catK) \neq \emptyset$, then the pair $(p\Spc(\catK), p\supp)$ is a final object in the category $\catS^{cm}(\catK)$. That is, if it is nonempty, $(p\Spc(\catK), p\supp)$ is a commutative support data, and for any other commutative support data $(Y, \sigma)$ on $\catK$, there exists a unique continuous map $f: Y \to p\Spc(\catK)$ such that $\sigma(x) = f\inv(p\supp(x)).$ Explicitly, the map is defined by \[f(x) = \{a \in \catK \mid x \not\in \sigma(a)\},\] on objects and $f(x)$ a full subcategory. 
	\end{theorem}
	\begin{proof}
		This verification follows essentially the same as the verification of Theorem \ref{thm:thm6}. First, one verifies that $(p\Spc(\catK), p\supp)$ is a commutative support data, which follows essentially the same as Proposition \ref{prop:multsuppdata} using Definition \ref{def:partialprime} - we omit this for brevity. Then, one must verify analogues of Proposition \ref{prop:t0} and Lemma \ref{lem:equalityoffuncs} - the proofs of these follow verbatim, as they are purely set-theoretic and do not rely on deeper topological properties of the spectrum. 
		
		Now, let $(Y,\sigma)$ be a commutative support datum for $\catK$ and let $x \in Y$. First, by definition we have that $f(x) \in p\supp(a)$ if and only if $a \not\in f(x)$ if and only if $x \in \sigma(a)$, so $f: Y \to \Spc^{cm}(\catK)$ is indeed the map specified and is unique. It remains to verify that $f$ is well-defined, i.e. that $f(x)$ is a partial prime $\otimes$-ideal. First, $f(x)$ is full by assumption. We have $f(x)$ is closed under shifts, $a \in f(x)$ if and only if $x \not\in \sigma(a)$ if and only if $x \not\in \sigma(\Sigma a)$ if and only if $\Sigma(a) \in f(x)$. Assume $a, b\in \catK$ commute. Then $a \oplus b \in f(x)$ if and only if $x \not\in \sigma(a \oplus b)$ if and only if $x \not \in \sigma(a) \cup \sigma(b)$ if and only if $a \in f(x)$ and $b \in f(x)$. Similarly, $a \otimes b \in f(x)$ if and only if $x \not\in \sigma(a \otimes b)$ if and only if $x \not\in \sigma(a) \cap \sigma(b)$ if and only if $a \in f(x)$ or $b \in f(x)$. Finally, let $a \to b\to c \to \Sigma a$ be a commuting triangle. If $a, b \in f(x)$ then $ x\not\in \sigma(a)$ and $x\not\in \sigma(b)$, hence $x\not\in \sigma(a) \cup \sigma(b)$. Therefore, we have (after shifting) $x\not\in \sigma(c)$, hence $c\in f(x)$ as well, as desired.  
	\end{proof}
	
	\begin{remark}
		In fact, any functorial spectral theory $(F,\sigma)$ on essentially small monoidal triangulated categories agreeing with $(\Spc,\supp)$ on all braided monoidal triangulated categories must also be a commutative support datum on $\catK$. Indeed, for any braided monoidal subcategory $\catC \subseteq \catK$ with the same tensor unit, one has a map on spectra $\Spc(\iota): F(\catK) \to F(\catC) \cong \Spc(\catC)$ satisfying \[ \Spc(\iota)\inv(\supp(x)) = \sigma(\iota(x)) = \sigma(x)\] for all $x \in \catC$.
		
		The next pertinent question is now: does there exist a monoidal triangulated category $\catK$ for which $p\Spc(\catK) = \emptyset$? At present, we are not aware of such a category, but conjecture such a category exists. We also ask if there exists a non-braided monoidal triangulated category $\catK$ for which $p\Spc(\catK) \cong \Spc\cp(\catK) \neq \emptyset$. Of course, if $p\Spc(\catK) = \emptyset$, then $\Spc\cp(\catK) = \emptyset$. Explicitly computing $p\Spc$ for sufficiently noncommutative monoidal triangulated categories appears challenging, since it remains unclear when an exact triangle is commutative. 
	\end{remark}

	\bibliography{bib}
	\bibliographystyle{alpha}
	
\end{document}